\documentclass[12pt, twoside, eqno]{article}
\usepackage{latexsym}
\usepackage{amssymb}
\usepackage{amsfonts}
\usepackage{amsmath}
\begin{document}\begin{center}

\noindent{\Large\bf From ordered semigroups to ordered 
hypersemigroups}\bigskip

\noindent{\bf Niovi KEHAYOPULU$^*$}

{\small Replacement to arXiv:1610.03880}
\bigskip

\noindent {\small \it Department of Mathematics,
University of Athens,
15784 Panepistimiopolis, Greece}\end{center}
{\small
\footnotetext {
\noindent *Correspondence: nkehayop@math.uoa.gr

\hspace{0.5cm}2010 {\it AMS Mathematics Subject Classification:} 
Primary 06F99, 20M99, 08A72;

\hspace{0.5cm}Secondary 06F05 }\bigskip

\noindent{\bf Abstract:}  We wrote this paper as an example to show 
the way we pass from ordered semigroups to ordered hypersemigroups. 
\medskip

  }
\section{Introduction and prerequisites}For an ordered semigroup 
$(S,\cdot,\le)$ and a subset $A$ of $S$, we denote by $(A]$ the 
subset of $S$ defined by $(A]:=\{t\in S \mid t\le a \mbox { for some 
} a\in A\}$. If $(S,\cdot,\le)$ is an ordered semigroup, a nonempty 
subset $A$ of $S$ is called a right (resp. left) ideal of $S$ if (1) 
$AS\subseteq A$ (resp. $SA\subseteq A$) and (2) if $a\in A$ and $S\ni 
b\le a$, then $b\in A$, that is, if $(A]=A$ (cf. e.g. [3]). An 
ordered semigroup $(S,\cdot,\le)$ is called regular [2] if for every 
$a\in S$ there exists $x\in S$ such that $a\le axa$. It is called 
intra-regular [3] if for every $a\in S$ there exist $x,y\in S$ such 
that $a\le xa^2 y$. A subset $B$ of an ordered semigroup 
$(S,\cdot,\le)$ is called a bi-ideal of $S$ [1] if (1) $BSB\subseteq 
B$ and (2) if $a\in B$ and $S\ni b\le a$, then $b\in B$. A subset $Q$ 
of $S$ is called a quasi-ideal of $S$ [4] if (1) $(QS]\cap 
(SQ]\subseteq Q$ and (2) if $a\in Q$ and $S\ni b\le a$, then $b\in 
Q$. We have seen in [6], that an ordered semigroup $S$ is regular if 
and only if for every right ideal $A$ and every left ideal $B$ of 
$S$, we have $A\cap B=(AB]$, equivalently, $A\cap B\subseteq (AB]$. 
In [8] we have shown that an ordered semigroup $S$ is intra-regular 
if and only if for every right ideal $A$ and every left ideal $B$ of 
$S$ we have $A\cap B\subseteq (BA]$. In [5], we proved that if $S$ is 
a regular ordered semigroup, then $B$ is a bi-ideal of $S$ if and 
only if there exist a right ideal $R$ and a left ideal $L$ of $S$ 
such that $B=(RL]$. We have seen in [6] that an ordered semigroup $S$ 
is regular if and only if the ideals of $S$ are idempotent and for 
every right ideal $A$ and every left ideal $B$ of $S$, the set $(AB]$ 
is a quasi-ideal of $S$. In [7] we studied some of the above results 
for hypersemigroups. In the present paper we examine the results on 
ordered semigroups mentioned above for ordered hypersemigroups. The 
results of ordered semigroups can be transferred to ordered 
hypersemigroups in the way indicated in the present paper. For the 
definitions and notations not given in the present paper we refer to 
[7]. The concept of an ordered semigroup can be transferred in a 
natural way to ordered hypersemigroups, for this concept we will 
refer to [9], and we will give further information about it in a next 
paper. In the proofs we tried to use sets instead of elements to show 
their pointless character.
\section{Main results} {\bf 1.} Based on ordered semigroups, we first 
introduce the concepts of regular and intra-regular hypersemigroups 
and we characterize them in terms of right and left ideals. The 
concept of regularity of ordered semigroups is naturally transferred 
to an hypersemigroup $H$ in the following definition\medskip

\noindent{\bf Definition 1.} An ordered hypersemigroup 
$(H,\circ,\le)$ is called {\it regular} if for every $a\in H$ there 
exist $x,y\in H$ such that $y\in (a\circ x)*\{a\}$ and $a\le y$.
\medskip

\noindent{\bf Proposition 2.} {\it Let H be an ordered 
hypersemigroup. The following are equivalent:

$(1)$ H is regular.

$(2)$ $a\in (\{a\}*H*\{a\}]$ for every $a\in H$.

$(3)$ $A\subseteq (A*H*A]$ for every $A\in {\cal P}^*(H)$}.\medskip

\noindent{\bf Definition 3.} Let $H$ be an ordered hypergroupoid. A 
nonempty subset $A$ of $H$ is called a {\it right} (resp. {\it left}) 
{\it ideal} of $H$ if

$(1)$ $A*H\subseteq A$ (resp. $H*A\subseteq A$) and

$(2)$ if $a\in A$ and $H\ni b\le a$, then $b\in A$, that is, if 
$(A]=A$. \medskip

\noindent{\bf Lemma 4.} {\it Let H be an ordered hypergroupoid and 
$A,B\in {\cal P}^*(H)$. Then we have$$(A]*(B]\subseteq (A*B].$$}
\noindent{\bf Lemma 5.} {\it For an hypersemigroup H and a nonempty 
subset A of H, we denote by $R(A)$ (resp. $L(A)$), the right (resp. 
left) ideal of H generated by A, and by $I(A)$ the ideal of $H$ 
generated by A. We have $$R(A)={\Big(}A\cup (A*H){\Big]}, \; 
L(A)={\Big(}A\cup (H*A){\Big]}, $$ $$I(A)={\Big(}A\cup (H*A)\cup 
(A*H)\cup (H*A*H){\Big]}.$$}{\bf Lemma 6.} {\it Let H be an 
hypergroupoid. If A is a right ideal and B a left ideal of H, then 
$A\cap B\not=\emptyset$.} \medskip

\noindent{\bf Lemma 7.} {\it Let H be an ordered hypergroupoid and 
$A,B\in {\cal P}^*(H)$. Then we 
have$${(A*B]=\Big(}(A]*(B]{\Big]}={\Big(}(A]*B{\Big]}= 
{\Big(}A*(B]{\Big]}.$$}{\bf Theorem 8.} {\it An ordered 
hypersemigroup H is regular if and only if for every right ideal A 
and every left ideal B of H, we have$$A\cap B=(A*B], \mbox { 
equivalently, } A\cap B\subseteq (A*B]).$$}{\bf Proof.} Let $H$ be 
regular, $A$ a right ideal and $B$ a left ideal of $H$. Since $A\cap 
B\not=\emptyset$ and $H$ is regular, by Proposition 2, we 
have\begin{eqnarray*}A\cap B&\subseteq&{\Big(}(A\cap B)*H*(A\cap 
B){\Big]}\subseteq {\Big(}(A*H)*B){\Big]}\subseteq (A*B] 
\\&\subseteq& (A*H]\cap (H*B]\subseteq (A]\cap (B]=A\cap 
B.\end{eqnarray*}Thus we have $A\cap B=(A*B]$.

Suppose now that $A\cap B\subseteq (A*B]$ for every right ideal $A$ 
and every left ideal $B$ of $H$, and let $A\in {\cal P}^*(H)$. Then 
we have\begin{eqnarray*}A&\subseteq&R(A)\cap 
L(A)\subseteq{\Big(}R(A)*L(A){\Big]}={\Bigg(}{\Big(}A\cup 
(A*H){\Big]}*{\Big(}A\cup 
(H*A){\Big]}{\Bigg]}\\&\subseteq&{\Bigg(}{\Big(}A\cup 
(A*H){\Big)}*{\Big(}A\cup (H*A){\Big)}{\Bigg]} \mbox { (by Lemma 
7)}\\&=&{\Big(}(A*A)\cup (A*H*A)\cup 
(A*H*H*A){\Big]}\\&=&{\Big(}(A*A)\cup 
(A*H*A){\Big]}.\end{eqnarray*}Then \begin{eqnarray*}A*A&\subseteq&
{\Big(}(A*A)\cup (A*H*A){\Big]}*(A]\subseteq {\Big(}(A*A*A)\cup 
(A*H*A*A){\Big]}\\&\subseteq&(A*H*A].\end{eqnarray*}
Then we have $$A\subseteq {\Big(}(A*H*A]\cup 
(A*H*A){\Big]}={\Big(}(A*H*A]{\Big]}=(A*H*A]$$and, by Proposition 2, 
$H$ is regular. $\hfill\Box$

The concept of intra-regularity of ordered semigroups is naturally 
transferred to an hypersemigroup $H$ in the following 
definition.\medskip

\noindent{\bf Definition 9.} An ordered hypersemigroup $H$ is called 
{\it intra-regular} if for every $a\in H$ there exist $x,y,t\in H$ 
such that $$t\in (x\circ a)*(a\circ y) \mbox { and } a\le 
t.$$\noindent{\bf Proposition 10.} {\it Let H be an ordered 
hypersemigroup. The following are equivalent:

$(1)$ H is intra-regular.

$(2)$ $a\in (H*\{a\}*\{a\}*H]$ for every $a\in H$.

$(3)$ $A\subseteq (H*A*A*H]$ for every $A\in {\cal 
P}^*(H)$}.\medskip

\noindent{\bf Theorem 11.} {\it An ordered hypersemigroup H is 
intra-regular if and only if for every right ideal A and every left 
ideal B of H, we have$$A\cap B\subseteq (B*A].$$}{\bf Proof.} 
$\Longrightarrow$. Let $A$ be a right ideal and $B$ a left ideal of 
$H$. By Lemma 6, $A\cap B\in {\cal P}^*(H)$. Since $H$ is 
intra-regular, by Proposition 10, we have\begin{eqnarray*}A\cap 
B&\subseteq& {\Big(}H*(A\cap B)*(A\cap B)*H{\Big]}\\&\subseteq 
&{\Big(}(H*B)*(A*H){\Big]} 
\\&\subseteq&B*A.\end{eqnarray*}$\Longleftarrow$. Just for a change, 
let us prove it using elements. The same can be proved using sets, as 
we did in Theorem 8. Let $a\in H$. By hypothesis, we 
have\begin{eqnarray*}a&\in&R(a)\cap L(a)\ 
\subseteq{\Big(}L(a)*R(a){\Big]}\\&=&{\Bigg (}{\Big(}\{a\}\cup 
(H*\{a\}){\Big]}*{\Big(}\{a\}\cup (\{a\}*H){\Big]}{\Bigg 
]}\\&\subseteq&{\Bigg(}{\Big(}\{a\}\cup 
(H*\{a\}){\Big)}*{\Big(}\{a\}\cup (\{a\}*H){\Big)}{\Bigg]} \mbox { 
(by Lemma 7)}\\&=&{\Big(}(a\circ a)\cup (H*\{a\}*\{a\})\cup 
(\{a\}*\{a\}*H)\cup (H*\{a\}*\{a\}*H){\Big]}. \end{eqnarray*}Then 
$a\in (a\circ a]$ or $a\in (H*\{a\}*\{a\}]$ or $a\in (\{a\}*\{a\}*H]$ 
or $a\in (H*\{a\}*\{a\}*H]$. If $a\in (a\circ a]$, 
then\begin{eqnarray*}a&\in& {\Big(}(a\circ a]*(a\circ 
a]{\Big]}\subseteq {\Big(}(a\circ a)*(a\circ a){\Big]} \mbox { (by 
Lemma 7)}\\&\subseteq& (H*\{a\}*\{a\}*H]\end{eqnarray*} then, by 
Proposition 10, $H$ is intra-regular. If $a\in (H*\{a\}*\{a\}]$, 
then\begin{eqnarray*}a&\in& 
{\Big(}H*(H*\{a\}*\{a\}]*\{a\}{\Big]}=(H*H*\{a\}*\{a\}*\{a\}]\\&\subseteq& 
(H*\{a\}*\{a\}*H].\end{eqnarray*}If $a\in (H*\{a\}*\{a\}]$ or $a\in 
(H*\{a\}*\{a\}*H]$ in a similar way, we have $a\in 
(H*\{a\}*\{a\}*H]$. In any case, $H$ is intra-regular. $\hfill\Box$

Below we first transfer the concept of bi-ideals from ordered 
semigroups to ordered hypersemigroups, then we characterize the 
bi-ideals of ordered hypersemigroups in terms of right and left 
ideals. \medskip

\noindent{\bf Definition 12.} Let $(H,\circ,\le)$ be an ordered 
hypersemigroup. An element\\$B\in P^* (H)$  is called a  {\it 
bi-ideal} of $H$ if

(1) $B*H*B\subseteq B$ and

(2) if $a\in B$ and $H\ni b\le a$, then $b\in B$, that is if 
$(B]=B$.\medskip

\noindent{\bf Proposition 13.} {\it Let H be an ordered 
hypersemigroup. If C is a right ideal and $D\in{\cal P}^*(H)$, then 
the set $B=(C*D]$ is a bi-ideal of H}.\medskip

\noindent{\bf Proof.} We 
have\begin{eqnarray*}B*H*B&=&(C*D]*(H]*(C*D]
\\&\subseteq& {\Big(}(C*D)*H{\Big]}*(C*D] \mbox { (by Lemma 
4)}\\&\subseteq&{\Big(}(C*D)*H*(C*D){\Big]} \mbox { (by Lemma 
7)}\\&=&{\Big(}(C*(D*H*C)*D{\Big]}\subseteq 
{\Big(}(C*H)*D{\Big]}\\&\subseteq& 
(C*D]=B.\end{eqnarray*}$\hfill\Box$\\In a similar way we have the 
following proposition\medskip

\noindent{\bf Proposition 14.} {\it Let H be an ordered 
hypersemigroup. If C is a left ideal and $D\in{\cal P}^*(H)$, then 
the set $B=(C*D]$ is a bi-ideal of H}.\\By Propositions 13 and 14 we 
have the following\medskip

\noindent{\bf Proposition 15.} {\it Let H be an ordered 
hypersemigroup. If C is a right ideal and D a left ideal of H, then 
the set $B=(C*D]$ is a bi-ideal of H}.\medskip

\noindent{\bf Proposition 16.} {\it Let H be a regular ordered 
hypersemigroup and B a bi-ideal of H. Then there exist a right ideal 
C and a left ideal D of H such that $B=(C*D]$.}\smallskip

\noindent{\bf Proof.} Since $B$ is a bi-ideal of $H$, we have 
$B*H*B\subseteq B$. Since $H$ is regular, by Proposition 2, we have 
$B\subseteq (B*H*B]$. Thus we have $B=(B*H*B]$.\\On the other 
hand,\begin{eqnarray*}R(B)*L(B)&=&{\Big(}(B\cup 
(B*H){\Big]}*{\Big(}(B\cup 
(H*B){\Big]}\\&\subseteq&{\Bigg(}{\Big(}B\cup 
(B*H){\Big)}*{\Big(}B\cup (H*B){\Big)}{\Bigg]} \mbox { (by Lemma 
7)}\\&=&{\Big(}(B*B)\cup (B*H*B)\cup 
(B*H*H*B){\Big]}\\&=&{\Big(}(B*B)\cup 
(B*H*B){\Big]}.\end{eqnarray*}Since $$B*B=(B]*(B*H*B]\subseteq 
(B*B*H*B]\subseteq (B*H*B],$$we 
have\begin{eqnarray*}R(B)*L(B)&\subseteq&{\Big(}(B*H*B]\cup 
(B*H*B){\Big]}={\Big(}(B*H*B]{\Big]}\\&=&(B*H*B]=B,\end{eqnarray*}
so ${\Big(}R(B)*L(B){\Big]}\subseteq (B]=B$.
On the other hand, we have \begin{eqnarray*}B&=&(B*H*B]\subseteq 
{\Bigg(}{\Big(}R(B)*H{\bigg)}*L(B){\bigg]}\\&\subseteq& 
{\Big(}R(B)*L(B){\Big]} \mbox { (since } R(B) \mbox { is a right 
ideal of } H).\end{eqnarray*}Hence we obtain 
$B={\Big(}R(B)*L(B){\Big]}$, where $R(B)$ is a right ideal and $L(B)$ 
is a left ideal of $H$. $\hfill\Box$\\By Propositions 15 and 16 we 
have the following theorem\smallskip

\noindent{\bf Theorem 17.} {\it Let H be a regular ordered 
hypersemigroup. An element\\$B\in {\cal P}^*(H)$ is a bi-ideal of H 
if and only if there exist a right ideal C and a left ideal D of H 
such that $B=(C*D]$.}\medskip

We characterize now the regular hypersemigroups in terms of right 
ideals, left ideals and quasi-ideals.
The concept of quasi-ideals of ordered semigroups is naturally 
transferred to ordered hypersemigroups as follows: \medskip

\noindent{\bf Definition 18.} Let $(H,\circ,\le)$ be an ordered 
hypersemigroup. An element\\$Q\in P^* (H)$  is called a  {\it 
quasi-ideal} of $H$ if

1) $(Q*H]\cap (H*Q]\subseteq Q$ and

2) if $a\in Q$ and $H\ni b\le a$, then $b\in Q$.\\
We first give the following definition, proposition.\medskip

\noindent{\bf Definition 19.} Let $H$ be an ordered hypersemigroup. 
An element $A\in {\cal P}^*(H)$ is called {\it subidempotent} if 
$(A*A]\subseteq A$. It is called {\it idempotent} if 
$(A*A]=A$.\medskip

\noindent{\bf Proposition 20.} {\it If H is an ordered hypergroupoid, 
then the right ideals and the left ideals of H are subidempotent. In 
particular, if H is a regular ordered hypersemigroup, then the right 
ideals and the left ideals of H are idempotent}.\medskip

\noindent{\bf Theorem 21.} {\it An ordered hypersemigroup H is 
regular if and only if the right ideals and the left ideals of H are 
idempotent and for every right ideal A and every left ideal B of H, 
the set $(A*B]$ is a quasi-ideal of H}.\medskip

\noindent{\bf Proof.} $\Longrightarrow$. By Proposition 20, the right 
and the left ideals of $H$ are idempotent. Let now $A$ be a right 
ideal and $B$ a left ideal of $H$. Since $H$ is regular, by Theorem 
8, we have $A\cap B=(A*B]$. It is  enough to prove that $A\cap B$ is 
a quasi-ideal of $H$. First of all, by Lemma 6, $A\cap B\in {\cal 
P}^*(H)$. Moreover,\begin{eqnarray*}{\Big(}(A\cap B)*H{\Big]}\cap 
{\Big(}H*(A\cap B){\Big]}&\subseteq&(A*H]\cap (H*B] 
\\&\subseteq&(A]\cap (B]=A\cap B,\end{eqnarray*}and if $x\in A\cap B$ 
and $H\ni y\le x$ then, since $y\le x\in A$ and $A$ is an ideal of 
$H$ we have $y\in A$ and since $y\le x\in B$, we have $y\in B$, so 
$y\in A\cap B$. Thus $A\cap B$ is a quasi-ideal of $H$, and so is 
$(A*B]$.\smallskip

\noindent$\Longleftarrow$. Let $A\in {\cal P}^*(H)$. Since $R(A)$ is 
a right ideal of $H$, it is idempotent, and we have 
have\begin{eqnarray*}A&\subseteq&R(A)={\Big(}R(A)*R(A){\Big]}=
{\Bigg(}{\Big(}A\cup (A*H){\Big]}*{\Big(}A\cup 
(A*H){\Big]}{\Bigg]}\\&=&{\Bigg(}{\Big(}A\cup 
(A*H){\Big)}*{\Big(}A\cup (A*H){\Big)}{\Bigg]} \mbox { (by Lemma 
7)}\\&=&{\Big(}(A*A)\cup (A*H*A)\cup (A*A*H)\cup (A*H*A*H){\Big]} 
\\&\subseteq&(A*H].\end{eqnarray*} Since $L(A)$ is a left ideal of 
$H$, in a similar way, we have $A\subseteq (H*A]$. Thus we get 
$A\subseteq (A*H]\cap (H*A]$. Since $(A*H]$ is a right ideal and 
$(H*A]$ a left ideal of $H$, they are idempotent, and we 
have\begin{eqnarray*}(A*H]\cap (H*A]&=&{\Big(}(A*H]*(A*H]{\Big]}\cap 
{\Big(}(H*A]*(H*A]{\Big]}\\
&=&{\Big(}(A*H)*(A*H){\Big]}\cap {\Big(}(H*A)*(H*A){\Big]} \mbox { 
(by Lemma 7)}\\&=&{\Big(}(A*H*A)*H{\Big]}\cap 
{\Big(}H*(A*H*A){\Big]}\\&=&{\Big(}(A*H*A]*H{\Big]}\cap 
{\Big(}H*(A*H*A]{\Big]} \mbox { (by Lemma 7)}. \end{eqnarray*}
On the other hand,\begin{eqnarray*}(A*H*A]&=&{\Big(}A*(H*H]*A{\Big]} 
\mbox { (since } H \mbox { is a left ideal of } H)\\
&=&{\Big(}A*(H*H)*A{\Big]} \mbox { (by Lemma 
7)}\\&=&{\Big(}(A*H)*(H*A){\Big]}\\&=&
{\Big(}(A*H]*(H*A]{\Big]} \mbox { (by Lemma 7)}.\end{eqnarray*}
Since $(A*H]$ is a right ideal and $(H*A]$ a left ideal of $H$, by 
hypothesis, the set ${\Big(}(A*H]*(H*A]{\Big]}$ is a quasi-ideal of 
$H$, and so is $(A*H*A]$. Hence we obtain 
$${\Big(}(A*H*A]*H{\Big]}\cap {\Big(}H*(A*H*A]{\Big]}\subseteq 
(A*H*A],$$then $A\subseteq (A*H*A]$ and by Proposition 2, $H$ is 
regular. $\hfill\Box$\medskip

\noindent {\bf 2.} We examine the results in [3] in case of ordered 
hypersemigroups. The well known concepts of prime and weakly prime 
and semiprime subsets of semigroups can be naturally transferred to 
hypergroupoids as follows:\medskip

\noindent{\bf Definition 22.} Let $H$ be an hypergroupoid (or an 
ordered hypergroupoid). A nonempty subset $T$ of $H$ is called a {\it 
prime subset} of $H$ if$$A, B\in {\cal P}^*(H) \mbox { and } 
A*B\subseteq T \mbox { implies } A\subseteq T \mbox { or } B\subseteq 
T.$${\bf Proposition 23.} {\it Let $(H,\circ)$ be an hypergroupoid. A 
nonempty subset T of H is a prime subset of H if and only if$$a,b\in 
H \mbox { such that } a\circ b\subseteq T \mbox { implies } a\in T 
\mbox { or } b\in T.$$}{\bf Proof.} $\Longrightarrow.$ Let $a,b\in 
H$, $a\circ b\subseteq T$. Since $\{a\}, \{b\}\in {\cal P}^*(H)$, 
$\{a\}*\{b\}=a\circ b\subseteq T$ and $T$ is prime, we have 
$\{a\}\subseteq T$ of $\{b\}\subseteq T$. Then $a\in T$ or $b\in 
T$.\\$\Longleftarrow$. Let $A, B\in {\cal P}^*(H)$ such that 
$A*B\subseteq T$ and let $A\not\subseteq T$ and $b\in B$. Take an 
element $a\in A$ such that $a\not\in T$. Since $a\circ b\subseteq 
A*B\subseteq T$, by hypothesis, we have $a\in T$ or $b\in T$. Since 
$a\not\in T$, we have $b\in T$. $\hfill\Box$\medskip

\noindent{\bf Definition 24.} Let $H$ be an hypergroupoid (or an 
ordered hypergroupoid). A nonempty subset $T$ of $H$ is called {\it 
weakly prime} if the following assertion is satified: $$\mbox { if } 
A, B \mbox { are ideals of } H \mbox { and } A*B\subseteq T, \mbox { 
then } A\subseteq T \mbox { or } B\subseteq T.$$
Recall that if $A$ and $B$ are ideals of $H$, then the sets $A$ and 
$B$ are nonempty, so $A*B\not=\emptyset$. \medskip

\noindent{\bf Definition 25.} Let $H$ be an hypergroupoid (or an 
ordered hypergroupoid). A nonempty subset $T$ of $H$ is called {\it 
semiprime} if$$\mbox { for any } A\in {\cal P}^*(H) \mbox { such that 
} A*A\subseteq T, \mbox { we have } A\subseteq T.$${\bf Proposition 
26.} {\it Let $H$ be an hypergroupoid (or an ordered hypergroupoid). 
A nonempty subset T of H is a semiprime subset of H if and only if 
$$\mbox { for any } a\in H \mbox { such that } a\circ a\subseteq T, 
\mbox { we have } a\in T.$$}{\bf Proof.} $\Longrightarrow.$ Let $a\in 
H$, $a\circ a\subseteq T$. Since $\{a\}\in {\cal P}^*(H)$ and 
$\{a\}*\{a\}=a\circ a\subseteq T$, by hypothesis, we have 
$\{a\}\subseteq T$, then $a\in T$.\\$\Longleftarrow$. Let $A\in {\cal 
P}^*(H)$ such that $A*A\subseteq T$ and let $a\in A$. Since $a\circ 
a\subseteq A*A\subseteq T$, by hypothesis, we have $a\in T$. 
$\hfill\Box$\medskip

\noindent{\bf Proposition 27.} {\it Let $(H,
 \circ,\le)$ be an ordered hypergroupoid. If A and B are ideals of H, 
then the intersection $A\cap B$ is an ideal of H as well.}\medskip

\noindent{\bf Proof.} First of all, since $A$ is a right ideal and 
$B$ a left ideal of $H$, we have $A\cap B\not=\emptyset$. Indeed: 
Take an element $a\in A$ and an element $b\in B$ 
$(A,B\not=\emptyset)$. Then $a\circ b\subseteq A*H\subseteq A$ and 
$a\circ b\subseteq B*H\subseteq B$, so $a\circ b\subseteq A\cap B$. 
Since $a\circ b\in {\cal P}^*(H)$, the set $a\circ b$ is a nonempty 
set. Thus $A\cap B$ is a nonempty subset of $H$ as well. In 
addition,$$(A\cap B)*H\subseteq A*H\subseteq A \mbox { and } (A\cap 
B)*H\subseteq B*H\subseteq B,$$so $(A\cap B)*H\subseteq A\cap B$, and 
$A\cap B$ is a right ideal of $(H,\circ)$. Similarly,$$H*(A\cap 
B)\subseteq H*A\subseteq A \mbox { and } H*(A\cap B)\subseteq 
H*B\subseteq B,$$so $H*(A\cap B)\subseteq A\cap B$, and $A\cap B$ is 
a left ideal of $(H,\circ)$. Let now $x\in A\cap B$ and $H\ni y\le x$ 
then, since $y\le x\in A$ we have $y\in A$ and, since $y\le x\in B$, 
we have $y\in B$, so $y\in A\cap B$. Thus $A\cap B$ is an ideal of 
$(H,\circ,\le)$. $\hfill\Box$\medskip

\noindent{\bf Theorem 28.} {\it Let $(H,\circ,\le)$ be an ordered 
hypergroupoid. The ideals of $H$ are idempotent if and only if for 
any two ideals A and B of H, we have$$A\cap B=(A*B].$$}{\bf Proof.} 
$\Longrightarrow$. Let $A$, $B$ be ideals of $H$. By Proposition 27, 
$A\cap B$ is an ideal of $H$. By hypothesis, we 
have\begin{eqnarray*}A\cap B&=&{\Big(}(A\cap B)*(A\cap 
B){\Big]}\subseteq (A*B]\\&\subseteq& (A*H]\cap (H*B]\subseteq 
(A]\cap (B]=A\cap B.\end{eqnarray*}Thus we have $A\cap 
B=(A*B]$.\\$\Longleftarrow$. Let $A$ be an ideal of $H$. By 
hypothesis, we have $A=A\cap A=(A*A]$, so $A$ is idempotent. 
$\hfill\Box$\medskip

\noindent{\bf Definition 29.} An ordered hypersemigroup 
$(H,\circ,\le)$ is called {\it semisimple} if, for every $a\in H$, we 
have$$a\in (H*\{a\}*H*\{a\}*H].$$This is equivalent to saying that 
for every $a\in H$ there exist $x,y,z,t\in H$ such that $t\in (x\circ 
a)*(y\circ a)*\{z\}$ and $a\le t$. Instead of writing $(x\circ 
a)*(y\circ a)*\{z\}$, we can write $\{x\}*(a\circ y)*(a\circ z)$ or  
$\{x\}*\{a\}*\{y\}*\{a\}*\{z\}$. On the other hand, $a\in 
(H*\{a\}*H*\{a\}*H]$ for every $a\in H$ if and only if $A\subseteq 
(H*A*H*A*H]$ for every nonempty subset $A$ of $H$. \medskip

\noindent{\bf Theorem 30.} {\it An ordered hypersemigroup 
$(H,\circ,\le)$ is semisimple if and only if the ideals of $H$ are 
idempotent.}\medskip

\noindent{\bf Proof.} $\Longrightarrow$. Let $A$ be an ideal of $H$. 
Then $A=(A*A]$. In fact: Let $x\in A$. By hypothesis, we have $x\in 
(H*\{x\}*H*\{x\}*H].$ Then $x\le t$ for some $t\in {\Big 
(}H*\{x\}*H{\Big )}*{\Big (}\{x\}*H{\Big)}.$ Then$$t\in a\circ b\; 
\mbox { for some } a\in H*\{x\}*H, \; b\in \{x\}*H.$$Since $a\in 
H*\{x\}*H\subseteq (H*A)*H\subseteq A*H\subseteq A$ and $b\in 
\{x\}*H\subseteq A*H\subseteq A$, we have $a\circ b\subseteq A*A$.
Since $x\le t\in A*A$, we have $x\in (A*A]$. Let now $x\in (A*A]$. 
Then $x\le t$ for some $t\in A*A$. Since $t\in A*A$, we have $t\in 
a\circ b$ for some $a,b\in A$. Since $a,b\in A$ and $A$ is a 
subsemigroup of $H$, we have $a\circ b\subseteq A*A\subseteq A$. 
Since $x\le t\in A$ and $A$ is an ideal of $H$, we have $x\in A$. 
Thus the ideals of $H$ are idempotent.\\$\Longleftarrow$. Let $a\in 
H$. By hypothesis, we have $I(a)=(I(a)*I(a)]$. In the implication 
$4)\Rightarrow 5)$ of Lemma 2 in [3], we replace the multiplication 
``$.$" by ``$*$", the proof follows. $\hfill\Box$\medskip

\noindent{\bf Proposition 31.} {\it Let $(H,\circ,\le)$ be an ordered 
hypersemigroup. If A is a left ideal and B a right ideal of H, then 
the set $(A*B]$ is an ideal of H}.\medskip

\noindent{\bf Proof.} Since $A$ and $B$ are nonempty subsets of $H$, 
the set $A*B$ is also a nonempty subset of $H$ and so is $(A*B]$. In 
addition,\begin{eqnarray*}H*(A*B]&=&(H]*(A*B]\subseteq{\Big(}H*(A*B){\Big]}\\
&=&{\Big(}(H*A)*B{\Big]}\subseteq (A*B],\end{eqnarray*}and
\begin{eqnarray*}(A*B]*H&=&(A*B]*(H]\subseteq{\Big(}(A*B)*H{\Big]}\\
&=&{\Big(}A*(B*H){\Big]}\subseteq (A*B].\end{eqnarray*}Let now $x\in 
(A*B]$ and $H\ni y\le x$. We have $x\le u$ for some $u\in A*B$. Since 
$H\ni y\le u\in A*B$, we have $y\in (A*B]$. Thus $(A*B]$ is an ideal 
of $H$. $\hfill\Box$\medskip

\noindent{\bf Corollary 32.} {\it If H is an ordered hypersemigroup 
and A,  B ideals of H, then the set $(A*B]$ is an ideal of 
H}.\medskip

\noindent{\bf Theorem 33.} {\it Let $(H,\circ,\le)$ be an ordered 
hypergroupoid. The ideals of H are weakly prime if and only if they 
are idempotent and they form a chain.}\medskip

\noindent{\bf Proof.} $\Longrightarrow$. Let $A$ be an ideal of $H$. 
Then $A=(A*A]$. Indeed: By Corollaty 32, the set $(A*A]$ is an ideal 
of $H$. Since $A*A\subseteq (A*A]$ and $(A*A]$ is weakly prime, we 
have $A\subseteq (A*A]\subseteq (A*H]\subseteq (A]=A$, so $A=(A*A]$. 
Let now $A$, $B$ be ideals of $H$. Then $A\subseteq B$ or $B\subseteq 
A$. Indeed: By Corollary 32, the set $(A*B]$ is an ideal of $H$. 
Since $A*B\subseteq (A*B]$ and $(A*B]$ is weakly prime, we have 
$A\subseteq (A*B]\subseteq (H*B]\subseteq (B]=B$ or $B\subseteq 
(A*B]\subseteq (A*H]\subseteq (A]=A$.

\noindent $\Longleftarrow$. If $T$ is an ideal of $H$, then $T$ is 
weakly prime. Indeed: Let $A, B$ be ideals of $H$ such that 
$A*B\subseteq T$. Since the ideals of $H$ are idempotent, by Theorem  
28, we have $A\cap B=(A*B].$ If $A\subseteq B$, then we have $A=A\cap 
B=(A*B]\subseteq (T]=T.$ If $B\subseteq A$, then $B=A\cap 
B=(A*B]\subseteq T$, so $T$ is weakly prime. $\hfill\Box$\medskip

\noindent{\bf Theorem 34.} {\it Let H be an ordered hypersemigroup. 
The ideals of H are prime if and only if they form a chain and H is 
intra-regular}.\medskip

\noindent{\bf Proof.} $\Longrightarrow$. Since the ideals of $H$ are 
prime, they are weakly prime. Then, by Theorem 33, they form a chain. 
Let now $a\in H$. Since the ideals of $H$ are semiprime, and 
${\Big(}\{a\}*\{a\}{\Big)}*{\Big(}\{a\}*\{a\}{\Big)}\subseteq 
{\Big(}H*\{a\}*\{a\}*H{\Big]}$, where ${\Big(}H*\{a\}*\{a\}*H{\Big]}$ 
is an ideal of $H$, we have $\{a\}*\{a\}\subseteq 
{\Big(}H*\{a\}*\{a\}*H{\Big]}$, and $a\in \{a\}\subseteq 
{\Big(}H*\{a\}*\{a\}*H{\Big]}$, thus $H$ is intra-regular.\smallskip

\noindent$\Longleftarrow$. Since $H$ is intra-regular, the ideals of 
$H$ are semiprime. In fact: Let $T$ be an ideal of $H$ and $a\in H$ 
such that $a\circ a\subseteq T$. Then$$a\in 
{\Big(}H*\{a\}*\{a\}*H{\Big]}={\Big(}H*(a\circ a)*H{\Big]}\subseteq 
(H*T*H]\subseteq (T]=T,$$ so $a\in T$, and $T$ is semiprime.
Since the ideals of $H$ are semiprime, the following assertions are 
satisfied:

(1) $I(A)=(H*A*H]$ for every $A\in {\cal P}^*(H)$. In fact:\\We have 
$(A*A)*(A*A)\subseteq (H*A*H]$, where $(H*A*H]$ is an ideal of $H$. 
Since $(H*A*H]$ is semiprime, we have $A*A\subseteq (H*A*H]$, and 
$A\subseteq (H*A*H]$, so $I(A)\subseteq (H*A*H]$.
On the other hand,
$$(H*A*H]\subseteq {\Big(}A\cup (H*A)\cup (A*H)\cup (H*A*H){\Big 
]}=I(A),$$
and condition (1) is satisfied.\smallskip

(2) $I(x\circ y)=I(x)\cap I(y)$ for every $x,y\in H$. In fact: Let 
$x,y\in H$.\\Since $x\circ y\subseteq I(x)*H\subseteq I(x)$, we have 
$I(x\circ y)\subseteq I(x)$. Since $x\circ y\subseteq H*I(y)\subseteq 
I(y)$, we have $I(x\circ y)\subseteq I(y)$. Thus we get $I(x\circ 
y)\subseteq I(x)\cap I(y)$. Let now $t\in I(x)\cap I(y)$. By (1), we 
have $t\in {\Big(}H*\{x\}*H{\Big]}$ and $t\in 
{\Big(}H*\{y\}*H{\Big]}$. Since$$H*\{x\}*H=\bigcup\limits_{a,b \in H} 
{\Big(}\{ a\} *\{ x\} *\{b\}{\Big)}  \mbox { and } 
H*\{y\}*H=\bigcup\limits_{c,d \in H} {\Big(}{\{ c\} *\{ y\} 
*\{d\}}{\Big)},$$we have$$t\le u, \mbox { where } u\in 
\{a\}*\{x\}*\{b\} \mbox { for some } a,b\in H$$and$$t\le v, \mbox { 
where } v\in \{c\}*\{y\}*\{d\} \mbox { for some } c,d\in H$$Then we 
have$$t\circ t\preceq v\circ 
u=\{v\}*\{u\}\subseteq\{c\}*{\Big(}\{y\}*\{d\}*\{a\}*\{x\}{\Big)}*\{b\}.$$On 
the other hand,$$\{y\}*\{d\}*\{a\}*\{x\}\subseteq I(x\circ 
y).$$Indeed: We 
have\begin{eqnarray*}{\Big(}\{y\}*\{d\}*\{a\}*\{x\}{\Big)}*{\Big(}
\{y\}*\{d\}*\{a\}*\{x\}{\Big)}&\subseteq& 
{\Big(}H*\{x\}*\{y\}*H{\Big)}\\&\subseteq&{\Big(}H*\{x\}*\{y\}*H{\Big]}
\\&=&I{\Big(}\{x\}*\{y\}{\Big)} \mbox { (by (1)) 
}.\end{eqnarray*}Since $I{\Big(}\{x\}*\{y\}{\Big)}$ is semiprime, we 
have ${\Big(}\{y\}*\{d\}*\{a\}*\{x\}{\Big)}\subseteq 
I{\Big(}\{x\}*\{y\}{\Big)}$. Since $I{\Big(}\{x\}*\{y\}{\Big)}$ is an 
ideal of $H$, we 
have$$\{c\}*{\Big(}\{y\}*\{d\}*\{a\}*\{x\}{\Big)}*\{b\}\subseteq 
H*I{\Big(}\{x\}*\{y\}{\Big)}*H\subseteq 
I{\Big(}\{x\}*\{y\}{\Big)}.$$Then $t\circ t\preceq v\circ u\subseteq 
I{\Big(}\{x\}*\{y\}{\Big)}$. Again since $I{\Big(}\{x\}*\{y\}{\Big)}$ 
is semiprime, we have$$t\in I{\Big(}\{x\}*\{y\}{\Big)}=I(x\circ 
y),$$and condition (2) holds.

Let now $T$ be an ideal of $H$ and $a,b\in H$ such that $a\circ 
b\subseteq T$. By hypothesis, we have $I(a)\subseteq I(b)$ or 
$I(b)\subseteq I(a)$. If $I(a)\subseteq I(b)$, then$$a\in 
I(a)=I(a)\cap I(b)=I(a\circ b)\subseteq I(T)=T,$$so $a\in T$. If
$I(b)\subseteq I(a)$, then$$b\in I(b)=I(a)\cap I(b)=I(a\circ 
b)\subseteq T,$$so $b\in T$. Thus $T$ is prime. $\hfill\Box$\medskip

\noindent{\bf Remark 35.} Here we give some examples of ordered 
hypersemigroups in which the ideals are idempotent. We first 
introduce the concept of regularity in ordered hypersemigroups. An 
ordered hypersemigroup $H$ is said to be {\it regular} if $a\in 
(\{a\}*H*\{a\}]$ for every $a\in H$, equivalently if $A\subseteq 
(A*H*A]$ for every $A\in {\cal P}^*(H)$. An ordered hypersemigroup 
$(H,\circ,\le)$ is regular if and only if for every $a\in H$ there 
exist $x,t\in H$ such that $t\in (a\circ x)*\{a\}$ and $a\le t$. 
Instead of writing $(a\circ x)*\{a\}$ we can write $\{a\}*(x\circ a)$ 
or $\{a\}*\{x\}*\{a\}$. In a regular hypersemigroup $H$ the right 
ideals and the left ideals are idempotent. In fact, let $A$ be a 
right ideal of $H$. Since $H$ is regular, we have$$A\subseteq 
((A*H)*A]\subseteq (A*A]\subseteq (A*H]\subseteq A,$$ so $(A*A]=A$. 
If $A$ is a left ideal of $H$ then, since $H$ is regular, we 
have$$A\subseteq (A*(H*A)]\subseteq (A*A]\subseteq (H*A]\subseteq 
A,$$ so $(A*A]=A$.
An ordered hypersemigroup $H$ is called {\it left regular} if $a\in 
(H*\{a\}*\{a\}]$ for every $a\in H$, equivalently if $A\subseteq 
(H*A*A]$ for every $A\in {\cal P}^*(H)$. An ordered hypersemigroup 
$(H,\circ,\le)$ is left regular if and only if, for every $a\in H$, 
there exist $x,t\in H$ such that $t\in \{x\}*(a\circ a)$ and $a\le 
t$. Instead of writing $\{x\}*(a\circ a)$, we can write 
$\{a\}*(x\circ a)$ or $\{a\}*\{x\}*\{a\}$. The left regular ordered 
hypersemigroups are intra-regular. Indeed, let $A\in {\cal P}^*(H)$. 
Then we have\begin{eqnarray*}A&\subseteq&(H*A*A]\subseteq 
{\Big(}H*(H*A*A]*A{\Big]}\\&\subseteq&{\Big(}(H]*(H*A*A]*(A]{\Big]}\\&=
&{\Big(}H*(H*A*A)*A{\Big]}\mbox { (by Lemma 
7)}\\&=&{\Big(}(H*H)*A*A*A{\Big]}\\&\subseteq&(H*A*A*H],\end{eqnarray*} 
so $H$ is intra-regular. In intra-regular ordered hypersemigroups the 
ideals are idempotent. In fact: Let $H$ be an intra-regular ordered 
hypersemigroup and $A$ an ideal of $H$. Then we have $A\subseteq 
(H*A*A*H]$, and so\begin{eqnarray*}(A*A]&\subseteq& 
{\Big(}(H*A*A*H]*A{\Big]}={\Big(}(H*A*A*H]*(A]{\Big]}\\&=&
{\Big(}(H*A*A*H)*A{\Big]}={\Big(}(H*A)*(A*H)*A{\Big]}\\&\subseteq&
(A*A*A]\subseteq (H*A]\subseteq (A]=A\\&\subseteq& 
{\Big(}(H*A)*(A*H){\Big]}\subseteq (A*A],\end{eqnarray*}so $(A*A]=A$, 
and $A$ is idempotent. An ordered hypersemigroup $H$ is said to be 
{\it right regular} if $a\in (\{a\}*\{a\}*H]$ for every $a\in H$, 
equivalently if $A\subseteq (A*A*H]$ for every $A\in {\cal P}^*(H)$. 
An ordered hypersemigroup $H$ is right regular if and only if for 
every $a\in H$, there exist $x,t\in H$ such that $t\in (a\circ 
a)*\{x\}$ and $a\le t$. Instead of writing $(a\circ a)*\{x\}$ we can 
write $\{a\}*(a\circ x)$ or $\{a\}*\{a\}*\{x\}$. The right regular 
ordered hypersemigroups are also intra-regular. Thus, in left 
regular, right regular or intra-regular ordered hypersemigroups the 
ideals are idempotent.\medskip

\noindent{\bf Remark 36.} In a similar way as in [1], we can prove 
the following:\\
(1) An ideal $T$ of an ordered hypersemigroup $H$ is weakly prime if 
and only if for all ideals $A$,$B$ of $H$ such that $(A*B]\cap 
(B*A]\subseteq T$, we have $A\subseteq T$ or $B\subseteq T$.\\(2) An 
ideal of an ordered hypersemigroup is prime if and only if it is both 
semiprime and weakly prime. In commutative ordered hypersemigroups 
the prime and the weakly prime ideals coincide.\medskip

\noindent{\bf 3.}  The theory of hypersemigroups and the theory of 
fuzzy hypersemigroups are parallel to each other, in the following 
sense: An hypersemigroup $H$ is intra-regular, for example, if and 
only if $A\cap B\subseteq B*A$ for every right ideal $A$ and every 
left ideal $B$ of $H$. And an hypersemigroup $H$ is intra-regular if 
and only if $f\wedge g\preceq g\circ f$ for every fuzzy right ideal 
$f$ and every fuzzy left ideal $g$ of $H$. An hypersemigroup $H$ is 
left quasi-regular if and only if $A\cap B\subseteq A*B$ for every 
ideal $A$ and every nonempty subset $B$ of $H$. And an hypersemigroup 
$H$ is left quasi-regular if and only if $f\wedge g\preceq f\circ g$ 
for every fuzzy ideal $A$ and every fuzzy subset $g$ of $H$. Further 
interesting information concerning this parallelism will be given in 
another paper. 

Following Zadeh, any mapping $f : H\rightarrow [0,1] $ of an 
hypergroupoid $H$ into the closed interval $[0,1]$ of real numbers is 
called a {\it fuzzy subset} of $H$ (or a {\it fuzzy set} in $H$) and 
the mapping $f_A$ (the so called characteristic function of $A$) is 
the fuzzy subset of $H$ defined as follows$$f_A : H \rightarrow 
\{0,1\} \mid x \rightarrow f_A (x)=\left\{ \begin{array}{l}
1\,\,\,\,\,$if$\,\,\,\,x \in A\\
0\,\,\,\,$if$\,\,\,\,x \notin A.
\end{array} \right.$$\\For an element $a$ of $H$, we denote by $A_a$ 
the subset of $H\times H$ defined by$$A_a:=\{(y,z)\in H\times H \mid 
a\in y\circ z\}.$$For two fuzzy subsets $f$ and $g$ of $H$, we denote 
by $f\circ g$ the fuzzy subset of $H$ defined as follows
$$f \circ g: H \to [0,1]\,\,\,a \to \left\{ \begin{array}{l}
\bigvee\limits_{(y,z) \in {A_a}} {\min \{ f(y),g(z)\} \,\,\,\,\mbox { 
if }\,\,\,{A_a} \ne \emptyset } \\
\,\,\,\,0\,\,\,\,\,\,\mbox { if }\,\,\,\,{A_a} = \emptyset. 
\,\,\,\,\,\,\,\,\,
\end{array} \right.$$As no confusion is possible, we denote the 
operation between fuzzy subsets of $H$ and the hyperoperation on $H$ 
by the same symbol. Denote by $F(H)$ the set of all fuzzy subsets of 
$H$ and by ``$\preceq$" the order relation on $F(H)$ defined 
by$$f\preceq g \;\Longleftrightarrow\; f(x)\le g(x) \mbox { for every 
} x\in H.$$

For two fuzzy subsets $f$ and $g$ of an hypergroupoid $H$ we denote 
by $f\wedge g$ the fuzzy subset of $H$ defined as follows$$f\wedge g 
: H\rightarrow [0,1] \mid x\rightarrow (f\wedge 
g)(x):=\min\{f(x),g(x)\}.$$One can easily see that the fuzzy subset 
$f\wedge g$ is the infimum of the fuzzy subsets $f$ and $g$, and this 
is why we write $f\wedge g=\inf\{f,g\}$. If $f$ is a fuzzy subset of 
$H$, then $f\wedge f=f$. Indeed, if $x\in H$, then $(f\wedge 
f)(x):=\min\{f(x),f(x)\}=f(x)$.

The concepts of fuzzy right and fuzzy left ideal of a semigroup due 
to Kuroki [12] can be naturally transferred to hypergroupoids as 
follows: A fuzzy subset $f$ of an hypergroupoid $H$ is called a {\it 
fuzzy right ideal} of $H$ if$$f(x\circ y)\ge f(x) \mbox { for every } 
x,y\in H,$$in the sense that if $x,y\in H$ and $u\in x\circ y$, then 
$f(u)\ge f(x)$.\\A fuzzy subset $f$ of an hypergroupoid $H$ is called 
a {\it fuzzy left ideal} of $H$ if$$f(x\circ y)\ge f(y) \mbox { for 
every } x,y\in H,$$meaning that if $x,y\in H$ and $u\in x\circ y$, 
then $f(u)\ge f(y)$.\\If $f$ is both a fuzzy right and a fuzzy left 
ideal of $H$, then it is called a {\it fuzzy ideal} of $H$. A fuzzy 
subset $f$ of an hypersemigroup $H$ is called a {\it fuzzy bi-ideal} 
of $H$ if$$f{\Big(}(x\circ y)*\{z\}{\Big)}\ge\min \{f(x),f(z)\} \mbox 
{ for every } x,y,z\in H,$$in the sense that if $x,y,z\in H$ and 
$u\in (x\circ y)*\{z\}$, then $f(u)\ge\min\{f(x),f(z)\}$.\\
If $H$ is an hypergroupoid, then $A$ is a right (resp. left) ideal of 
$H$ if and only if the characteristic function $f_A$ is a fuzzy left 
(resp. fuzzy right) ideal of $H$. If $H$ is an hypersemigroup then 
$A$ is a bi-ideal of $H$ if and only if $f_A$ is a fuzzy bi-ideal of 
$H$.\medskip

\noindent{\bf Theorem 37.} {\it Let $(H,\circ)$ be an hypersemigroup. 
The following are equivalent:\begin{enumerate}
\item H is regular.
\item $A\cap B=A*B$ for every right ideal A and every left ideal B of 
H.
\item $A\cap B\subseteq A*B$ for every right ideal A and every left 
ideal B of H.
\item $R(A)\cap L(A)\subseteq R(A)*L(A)$ for every $A\in
{\cal P}^*(H)$.
\item $R(a)\cap L(a)\subseteq R(a)*L(a)$ for every $a\in 
H$.\end{enumerate}}
\noindent{\bf Proof.} $(1)\Longrightarrow (2)$. Let $A$ be a right 
ideal and $B$ a left ideal of $H$. The set $A\cap B$ is a nonempty 
subset of $H$. Indeed: Take an element $a\in A$ and an element $b\in 
B$ $(A, B\not=\emptyset)$. Then $a\circ b\subseteq A*H\subseteq A$ 
and $a\circ b\subseteq H*B\subseteq B$, so $a\circ b\subseteq A\cap 
B$. Since $a\circ b\not=\emptyset$, we have $A\cap B\not=\emptyset$. 
Since $H$ is regular and $A\cap B\in {\cal P}^*(H)$, we 
have\begin{eqnarray*}A\cap B&\subseteq&(A\cap B)*H*(A\cap B)\subseteq 
A*H*B=(A*H)*B\\&\subseteq& A*B\subseteq (A*H)\cap (H*B)\subseteq 
A\cap B.\end{eqnarray*}Thus we have $A\cap B=A*B$.\\The implications 
$(2)\Rightarrow (3)\Rightarrow (4)\Rightarrow (5)$ are 
obvious.\\$(5)\Longrightarrow (1)$. Let $a\in H$. Since $R(a)$ is a 
right ideal of $H$ and $L(a)$ is a left ideal of $H$, by hypothesis, 
we have\begin{eqnarray*}a&\in&R(a)\cap L(a)\subseteq 
R(a)*L(a)={\bigg(}\{a\}\cup 
{\Big(}\{a\}*H{\Big)}{\bigg)}*{\bigg(}\{a\}\cup 
{\Big(}H*\{a\}{\Big)}{\bigg)}\\&=&(a\circ a)
\cup {\Big(}\{a\}*H*\{a\}{\Big)}\cup 
{\Big(}\{a\}*H*H*\{a\}{\Big)}\\&=&(a\circ a)\cup 
{\Big(}\{a\}*H*\{a\}{\Big)}.\end{eqnarray*}We have $a\in a\circ a$, 
so $a\in (a\circ a)*\{a\}$ or $a\in \{a\}*H*\{a\}$. In each case, $H$ 
is regular. $\hfill\Box$

It might be noted that in the above theorem, if $H$ is regular, then 
for every right ideal $A$ and every subset $B$ of $H$ or for every 
subset $A$ and every left ideal $B$ of $H$ we also have $A\cap 
B=A*B$. \medskip

\noindent{\bf Theorem 38.} {\it Let $(H,\circ)$ be an hypersemigroup. 
The following are equivalent:\begin{enumerate}
\item H is regular.
\item $f\wedge g=f\circ g$ for every fuzzy right ideal f and every 
fuzzy left ideal g of H.
\item $f\wedge g\preceq f\circ g$ for every fuzzy right ideal f and 
every fuzzy left ideal g of H.\end{enumerate}}
\noindent{\bf Proof.} $(1)\Longrightarrow (2)$. Let $f$ be a fuzzy 
right ideal and $g$ a fuzzy left ideal of $H$. Then we have $f\circ 
g\preceq f\wedge g$. In fact: Let $a\in H$. Then $(f\circ g)(a)\le 
(f\wedge g)(a)$. Indeed: If $A_a=\emptyset$, then $(f\circ g)(a):=0$. 
Since $a\in H$ and $f\wedge g$ is a fuzzy subset of $H$, we have 
$(f\wedge g)(a)\ge 0$, thus we have $(f\circ g)(a)\le (f\wedge 
g)(a)$. Let now $A_a\not=\emptyset$. 
Then\begin{equation}\tag{$\ast$}(f\circ g)(a):=\bigvee\limits_{(x,y) 
\in {A_a}} {\min \{ f(x),g(y)\} }\end{equation} We 
have\begin{equation}\tag{$\ast\ast$}\min\{f(x),g(y)\}\le (f\wedge 
g)(a) \mbox { for every } (x,y)\in A_a \end{equation}Indeed: Let 
$(x,y)\in A_a$. Then $a\in x\circ y$. Since $f$ is a fuzzy right 
ideal of $H$, we have $f(x\circ y)\ge f(x)$, then we have $f(a)\ge 
f(x)$. Since $g$ is a fuzzy left ideal of $H$, we have $g(x\circ 
y)\ge g(y)$, then $g(a)\ge g(y)$, so$$(f\wedge 
g)(a):=\min\{f(a),g(a)\}\ge \min\{f(x),g(y)\},$$ and condition $(**)$ 
is satisfied. By $(**)$, we have $$\bigvee\limits_{(x,y) \in {A_a}} 
\min \{ f(x),g(y)\}\le (f\wedge g)(a).$$Then, by $(*)$, $(f\circ 
g)(a)\le (f\wedge g)(a)$. Moreover, since $(H,\circ)$ is a regular 
hypersemigroup, $f$ a fuzzy right ideal of $H$ and $g$ a fuzzy subset 
of $H$, we have $f\wedge g\preceq f\circ g.$ In fact: Let $a\in H$. 
Since $H$ is regular, there exists $x\in H$ such that $a\in (a\circ 
x)*\{a\}$. Then $a\in u\circ a$ for some $u\in a\circ x$. Since $a\in 
u\circ a$, we have $(u,a)\in A_a$, then $$(f\circ 
g)(a):=\bigvee\limits_{(y,z) \in {A_a}} {\min \{ f(y),g(z)\} }\ge 
\min\{f(u),g(a)\}.$$Since $f$ is a fuzzy right ideal of $H$, we have 
$f(a\circ x)\ge f(a)$. Since $u\in a\circ x$, we have $f(u)\ge f(a)$. 
Then we have$$(f\circ g)(a)\ge 
\min\{f(u),g(a)\}\ge\min\{f(a),g(a)\}:=(f\wedge g)(a),$$so $f\wedge 
g\preceq f\circ g$. Thus condition (2) is satisfied.\\The implication 
$(2)\Rightarrow (3)$ is obvious.\\$(3)\Longrightarrow (1)$. 
By Theorem 37, it is enough to prove that $R(a)\cap L(a)\subseteq 
R(a)*L(a)$ for every $a\in H$. So, let $a\in H$ and $b\in R(a)\cap 
L(a)$. Since $R(a)$ is a right ideal of $H$, the characteristic 
function $f_{R(a)}$ is a fuzzy right ideal of $H$ and, since $L(a)$ 
is a left ideal of $H$,
$f_{L(a)}$ is a fuzzy left ideal of $H$. By hypothesis, we have 
$f_{R(a)}\wedge f_{L(a)}\preceq f_{L(a)}\circ f_{R(a)}$, then 
$${\Big(}f_{R(a)}\wedge f_{L(a)}{\Big)}(b)\le {\Big(}f_{R(a)}\circ 
f_{L(a)}{\Big)}(b).$$Thus$$\min\{f_{R(a)}(b), f_{L(a)}(b)\}\le 
{\Big(}f_{R(a)}\circ f_{L(a)}{\Big)}(b).$$
Since $b\in R(a)$ and $b\in L(a)$, we have 
$f_{R(a)}(b)=f_{L(a)}(b)=1$, and so $$1\le {\Big(}f_{R(a)}\circ 
f_{L(a)}{\Big)}(b).$$If $A_b=\emptyset$, then ${\Big(}f_{R(a)}\circ 
f_{L(a)}{\Big)}(b)=0$ which is impossible. Thus we have 
$A_b\not=\emptyset$ and$${\Big(}f_{R(a)}\circ 
f_{L(a)}{\Big)}(b)=\bigvee\limits_{(y,z) \in {A_b}} {\min \{ 
{f_{R(a)}}(y),{f_{L(a)}}(z)\} }.$$Then there exists $(y,z)\in A_b$ 
such that $y\in R(a)$ and $z\in L(a)\hfill(*)$\\Indeed: Suppose there 
is no $(y,z)\in A_b$ such that $y\in R(a)$ and $z\in L(a)$. Then, for 
each $(y,z)\in A_b$, we have $y\notin R(a)$ or $z\notin L(a)$. Then, 
for each $(y,z)\in A_b$, we have $f_{R(a)}(y)=0$ or $f_{L(a)}(z)=0$. 
Then $\min\{f_{R(a)}(y), f_{R(a)}(z)\}=0$ for every $(y,z)\in A_b$. 
Then
$\bigvee\limits_{(y,z) \in {A_b}} {\min \{ 
{f_{R(a)}}(y),{f_{L(a)}}(z)\} }=0$, so ${\Big(}f_{R(a)}\circ 
f_{L(a)}{\Big)}(b)=0$ which is no possible. By $(*)$, we have
$b\in y\circ z\subseteq R(a)*L(a)$, and the proof is complete. 
$\hfill\Box$\medskip

\noindent{\bf Theorem 39.} (cf. also [18]) {\it Let $(H,\circ)$ be an 
hypersemigroup. The following are equivalent:\begin{enumerate}
\item H is intra-regular.
\item $A\cap B\subseteq B*A$ for every right ideal A and every left 
ideal B of H.
\item $R(A)\cap L(A)\subseteq L(A)*R(A)$ for every $A\in
{\cal P}^*(H)$.
\item $R(a)\cap L(a)\subseteq L(a)*R(a)$ for every $a\in 
H$.\end{enumerate}}

\noindent{\bf Theorem 40.}  {\it An hypersemigroup $(H,\circ)$ is 
intra-regular if and only if, for every fuzzy right ideal f and every 
fuzzy left ideal g of H, we have $f\wedge g\preceq g\circ 
f.$}\medskip

\smallskip

\noindent{\bf Proof.} $\Longrightarrow$. Let $a\in H$. Since $H$ is 
intra-regular, there exist $x,y\in H$ such that $a\in (x\circ 
a)*(a\circ y)$. Then $a\in u\circ v$ for some $u\in x\circ a$, $v\in 
a\circ y$. Since $a\in u\circ v$, we have $(u,v)\in A_a$, then we 
have$$(g\circ f)(a):=\bigvee\limits_{(h,k) \in {A_a}} {\min \{ 
g(h),f(k)\}}\ge \min\{g(u),f(v)\}.$$Since $g$ is a fuzzy left ideal 
of $H$, we have $g(x\circ a)\ge g(a)$. Since $u\in x\circ a$, we get 
$g(u)\ge g(a)$. Since $f$ is a fuzzy right ideal of $H$, we have 
$f(a\circ y)\ge f(a)$. Since $v\in a\circ y$, we have $f(v)\ge f(a)$. 
Thus we have$$(g\circ f)(a)\ge \min\{g(u),f(v)\}\ge 
\min\{g(a),f(a)\}=(f\wedge g)(a),$$thus $f\wedge g\preceq g\circ 
f$.\\$\Longleftarrow$. By Theorem 39, it is enough to prove that 
$R(a)\cap L(a)\subseteq L(a)*R(a)$ for every $a\in H$. Let now $a\in 
H$ and $b\in R(a)\cap L(a)$. As $f_{R(a)}$ is a fuzzy right ideal and 
$f_{L(a)}$ is a fuzzy left ideal of $H$, by hypothesis, we have 
$${\Big(}f_{R(a)}\wedge f_{L(a)}{\Big)}(b)\le {\Big(}f_{L(a)}\circ 
f_{R(a)}{\Big)}(b).$$Thus$$\min\{f_{R(a)}(b), f_{L(a)}(b)\}\le 
{\Big(}f_{L(a)}\circ f_{R(a)}{\Big)}(b).$$Since $b\in R(a)$ and $b\in 
L(a)$, we have $f_{R(a)}(b)=f_{L(a)}(b)=1$, and so $$1\le 
{\Big(}f_{L(a)}\circ f_{R(a)}{\Big)}(b).$$If $A_b=\emptyset$, then 
${\Big(}f_{L(a)}\circ f_{R(a)}{\Big)}(b)=0$ which is impossible. Then 
we have $A_b\not=\emptyset$ and then $${\Big(}f_{L(a)}\circ 
f_{R(a)}{\Big)}(b)=\bigvee\limits_{(y,z) \in {A_b}} {\min \{ 
{f_{L(a)}}(y),{f_{R(a)}}(z)\} }.$$Then there exists $(y,z)\in A_b$ 
such that $y\in L(a)$ and $z\in R(a)$. Then we get $b\in y\circ 
z\subseteq L(a)*R(a)$. $\hfill\Box$

The concept of left quasi-regular semigroups can be naturally 
transferred to hypersemigroups in the definition below.\medskip

\noindent{\bf Definition 41.} An hypersemigroup $(H,\circ)$ is called 
{\it left quasi-regular} if for every $a\in H$ there exist $x,y\in H$ 
such that $a\in (x\circ a)*(y\circ a)$.\medskip

\noindent{\bf Proposition 42.} {\it Let $(H,\circ)$ be an 
hypersemigroup. The following are equivalent:\begin{enumerate}
\item H is left quasi-regular.
\item $a\in H*\{a\}*H*\{a\}$ for every $a\in H$.
\item $A\subseteq H*A*H*A$ for every $A\in {\cal 
P}^*(H)$.\end{enumerate}}
\noindent{\bf Theorem 43.} {\it Let $(H,\circ)$ be an hypersemigroup. 
The following are equivalent:\begin{enumerate}
\item H is left quasi-regular.
\item $A\cap B\subseteq A*B$ for every ideal $A$ and every nonempty 
subset $B$ of $H$.
\item $A\cap B\subseteq A*B$ for every ideal $A$ and every bi-ideal 
$B$ of $H$.
\item $A\cap B\subseteq A*B$ for every ideal $A$ and every left ideal 
$B$ of $H$.
\item $I(A)\cap L(A)\subseteq I(A)*L(A)$ for every $A\in
{\cal P}^*(H)$.
\item $I(a)\cap L(a)\subseteq I(a)*L(a)$ for every $a\in 
H$.\end{enumerate}}
\noindent{\bf Proof.} $(1)\Longrightarrow (2)$. Let $A$ be an ideal, 
$B$ a nonempty subset of $H$ and $a\in A\cap B$. Since $a\in H$ and 
$H$ is left quasi-regular, by Proposition 42$(1)\Rightarrow (2)$, we 
have
$$a\in H*\{a\}*H*\{a\}={\Big(}H*\{a\}*H{\Big)}*\{a\}\subseteq 
(H*A*H)*B.$$ We have $H*A*H=(H*A)*H\subseteq A*H$ since $A$ is a left 
ideal of $H$ and $A*H\subseteq A$ since $A$ is a right ideal of $H$. 
Thus we have $H*A*H\subseteq A$. Then $a\in (H*A*H)*B\subseteq 
A*B$.\\The implications $(2)\Rightarrow (3)$ and $(4)\Rightarrow 
(5)\Rightarrow (6)$ are obvious, and $(3)\Rightarrow (4)$ since the 
left ideals of $H$ are bi-ideals of $H$ as well.\\$(6)\Longrightarrow 
(1)$. Let $a\in H$. By hypothesis, we have\begin{eqnarray*}a&\in& 
I(a)\cap L(a)\subseteq
I(a)*L(a)\\&=&{\bigg(}\{a\}\cup {\Big(}H*\{a\}{\Big)}\cup 
{\Big(}\{a\}*H{\Big)}\cup 
{\Big(}H*\{a\}*H{\Big)}{\bigg)}*{\bigg(}\{a\}\cup 
{\Big(}H*\{a\}{\Big)}{\bigg)}\\&=&{\Big(}\{a\}*\{a\}{\Big)}\cup 
{\Big(}H*\{a\}*\{a\}{\Big)}\cup {\Big(}\{a\}*H*\{a\}{\Big)}\cup 
{\Big(}H*\{a\}*H*\{a\}{\Big)}. \end{eqnarray*}If $a\in \{a\}*\{a\}$, 
then \begin{eqnarray*}a\in \{a\}&\subseteq&\{a\}*\{a\}\subseteq 
{\Big(}\{a\}*\{a\}{\Big)}*{\Big(}\{a\}*\{a\}{\Big)}\\
&\subseteq& H*\{a\}*H*\{a\}.
\end{eqnarray*}If $a\in H*\{a\}*\{a\}$, then\begin{eqnarray*}a\in 
\{a\}&\subseteq& H*\{a\}*{\Big(}H*\{a\}*\{a\}{\Big)}\subseteq 
H*\{a\}*(H*H)*\{a\}\\
&\subseteq&H*\{a\}*H*\{a\}.\end{eqnarray*}If $a\in \{a\}*H*\{a\}$, 
then{\begin{eqnarray*}\{a\}&\subseteq&\{a\}*H*{\Big(}\{a\}*H*\{a\}{\Big)}\subseteq
(H*H)*{\Big(}\{a\}*H*\{a\}{\Big)}\\&\subseteq&H*\{a\}*H*\{a\}.\end{eqnarray*}
In each case, $a\in H*\{a\}*H*\{a\}$, so $H$ is left quasi-regular. 
$\hfill\Box$\medskip

\noindent{\bf Theorem 44.} {\it An hypersemigroup $(H,\circ)$ is left 
quasi-regular if and only if, for any left ideals $A$ and $B$ of $H$, 
we have $A\cap B\subseteq A*B$.}\medskip

\noindent{\bf Proof.} $\Longrightarrow$. Let $A$, $B$ be left ideals 
of $H$ and $a\in A\cap B$. Since $H$ is left quasi-regular, there 
exist $x,y\in H$ such that\begin{eqnarray*}a\in (x\circ a)*(y\circ 
a)&=&\{x\}*\{a\}*\{y\}*\{a\}\\&\subseteq&(H*A)*(H*B)\subseteq 
A*B.\end{eqnarray*}
$\Longleftarrow$. Let $A\in {\cal P}^*(H).$ Then $A\subseteq 
H*A*H*A$. In fact, by hypothesis, we have\begin{eqnarray*}
A&\subseteq&L(A)=L(A)\cap L(A)\subseteq L(A)*L(A)={\Big(}A\cup 
(H*A){\Big)}*{\Big(}A\cup (H*A){\Big)}\\&=&(A*A)\cup (H*A*A)\cup 
(A*H*A)\cup (H*A*H*A).\end{eqnarray*}Then we 
have\begin{eqnarray*}A*A&\subseteq&(A*A*A)\cup (A*H*A*A)\cup 
(A*A*H*A)\cup (A*H*A*H*A)\\&\subseteq&A*H*A,\end{eqnarray*}
from which $H*A*A\subseteq H*A*H*A$. Thus we obtain$$A\subseteq 
(A*H*A)\cup (H*A*H*A),$$ then we get\begin{eqnarray*}A*H*A&\subseteq& 
(A*H*A*H*A)\cup (H*A*H*A*H*A)\\&\subseteq& H*A*H*A,\end{eqnarray*}so 
$A\subseteq H*A*H*A$, and $H$ is left quasi-regular. $\hfill\Box$

A subset $A$ of an hypergroupoid $(H,\circ)$ is called {\it 
idempotent} if $A*A=A$.\medskip

\noindent{\bf Theorem 45.} {\it An hypersemigroup $(H,\circ)$ is left 
quasi-regular if and only if the left ideals of $H$ are 
idempotent.}\medskip

\noindent{\bf Proof.} $\Longrightarrow$. If $L$ be a left ideal of 
$H$ then, since $L\in {\cal P}^*(H)$ and $H$ is left quasi-regular, 
by Proposition $42(1)\Rightarrow (3)$, we have$$L\subseteq 
(H*L)*(H*L)\subseteq L*L\subseteq L.$$ The same is also a consequence 
of the $\Rightarrow$-part of the Theorem 44. Indeed, if $L$ is a left 
ideal of $H$ then, by Theorem 44, we have $L\subseteq L*L\subseteq 
H*L\subseteq L$, so $L*L=L$.\\$\Longleftarrow$. By Theorem 44, it is 
enough to prove that for every ideal $A$ and every left ideal $B$ of 
$H$, we have $A\cap B\subseteq A*B$. Let now $A$ be an ideal and $B$ 
a left ideal of $H$. The set $A\cap B$ is a nonempty subset of $H$ 
and $H*(A\cap B)\subseteq (H*A)\cap (H*B)\subseteq A\cap B$, so the 
set $A\cap B$ is a left ideal of $H$. By hypothesis, we have $A\cap 
B=(A\cap B)*(A\cap B)\subseteq A*B$. $\hfill\Box$\medskip

\noindent{\bf Theorem 46.} {\it Let $(H,\circ)$ be an hypersemigroup. 
The following are equivalent:\begin{enumerate}
\item H is left quasi-regular.
\item $f\wedge g\preceq f\circ g$ for every fuzzy ideal f and every 
fuzzy subset g of H.
\item $f\wedge g\preceq f\circ g$ for every fuzzy ideal f and every 
fuzzy bi-ideal g of H.
\item $f\wedge g\preceq f\circ g$ for every fuzzy ideal f and every 
fuzzy left ideal g of H.\end{enumerate}}
\noindent$(1)\Longrightarrow (2)$. Let $f$ be a fuzzy ideal, $g$ a 
fuzzy subset of $H$ and $a\in H$. Since $H$ is left quasi-regular, 
there exist $x,y\in H$ such that $a\in{\Big (}(x\circ a)*\{y\}{\Big 
)}*\{a\}$. Then $a\in u\circ a$ for some $u\in (x\circ a)*\{y\}$. In 
addition, $u\in v\circ y$ for some $v\in x\circ a$. On the other 
hand, since $(u,a)\in A_a$, we have$$(f\circ 
g)(a):=\bigvee\limits_{(h,k) \in {A_a}} {\min \{ f(h),g(k)\}}\ge \min 
\{f(u),g(a)\}.$$Since $f$ is a fuzzy right ideal of $H$, we have 
$f(v\circ y)\ge f(v)$ and since $u\in v\circ y$, we have $f(u)\ge 
f(v)$. Since $f$ is a fuzzy left ideal of $H$, we have $f(x\circ 
a)\ge f(a)$ and since $v\in x\circ a$, we have $f(v)\ge f(a)$. Thus 
we have $f(u)\ge f(a)$, and$$(f\circ g)(a)\ge 
\min\{f(a),g(a)\}=(f\wedge g)(a),$$so $f\wedge g\preceq f\circ 
g$.\\The implication $(2)\Rightarrow (3)$ is obvious and
$(3)\Rightarrow (4)$ since every fuzzy left ideal is a fuzzy bi-ideal 
of $H$.\\$(4)\Longrightarrow (1)$. By Theorem 43, it is enough to 
prove that $I(a)\cap L(a)\subseteq I(a)*L(a)$ for every $a\in H$. Let 
now $a\in H$ and $b\in I(a)\cap L(a)$. Since $I(a)$ is an ideal of 
$H$, the characteristic function $f_{I(a)}$ is a fuzzy ideal of $H$ 
and, since $L(a)$ is a left ideal of $H$, $f_{L(a)}$ is a fuzzy left 
ideal of $H$. By hypothesis, we have $f_{I(a)}\wedge f_{L(a)}\preceq 
f_{I(a)}\circ f_{L(a)}$, and so ${\Big(}f_{I(a)}\wedge 
f_{L(a)}{\Big)}(b)\le {\Big(}f_{I(a)}\circ f_{L(a)}{\Big)}(b)$, that 
is $$\min \{f_{I(a)}(b),f_{L(a)}(b)\}\le {\Big(}f_{I(a)}\circ 
f_{L(a)}{\Big)}(b).$$Since $b\in I(a)$, we have $f_{I(b)}(b)=1$ and, 
since $b\in L(a)$, we have $f_{L(b)}(b)=1$, so $\min 
\{f_{I(a)}(b),f_{L(a)}(b)\}=1$, and so $1\le {\Big(}f_{I(a)}\circ 
f_{L(a)}{\Big)}(b)$. If $A_b=\emptyset$, then ${\Big(}f_{I(a)}\circ 
f_{L(a)}{\Big)}(b)=0$ which is impossible. Thus we have 
$A_b\not=\emptyset$. Then $${\Big(}f_{I(a)}\circ 
f_{L(a)}{\Big)}(b)=\bigvee\limits_{(y,z) \in {A_b}} {\min 
\{f_{I(a)}(y),f_{L(a)}(z)\}}.$$Then there exists $(y,z)\in A_b \mbox 
{ such that } y\in I(a) \mbox { and } z\in L(a)\hfill(*)$\\Indeed: 
Suppose there is no $(y,z)\in A_b$ such that $y\in I(a)$ and $z\in 
L(a)$. Then, for every $(y,z)\in A_b$ we have $y\notin I(a)$ or 
$z\notin L(a)$. Then, for each $(y,z)\in A_b$ we have $f_{I(a)}(y)=0$ 
or $f_{L(a)}(z)=0$, so for each $(y,z)\in A_b$, we have $\min 
\{f_{I(a)}(y),f_{L(a)}(z)\}=0$, then ${\Big(}f_{I(a)}\circ 
f_{L(a)}{\Big)}(b)=0$ which is no possible.\\By $(*)$, we have $b\in 
y\circ z\subseteq I(a)*L(a)$. $\hfill\Box$\medskip

\noindent{\bf Theorem 47.} {\it An hypersemigroup $(H,\circ)$ is left 
quasi-regular if and only if for any fuzzy left ideals f and g of H, 
we have $f\wedge g\preceq f\circ g$.}\medskip

\noindent{\bf Proof.} $\Longrightarrow$. Let $f$ and $g$ be fuzzy 
left ideals of $H$ and $a\in H$. Then $(f\wedge g)(a)\le (f\circ 
g)(a)$. Indeed: By hypothesis, we have $a\in (x\circ a)*(y\circ a)$, 
so we have $a\in u\circ v$ for some $u\in x\circ a$, $v\in y\circ a$. 
Since $(u,v)\in A_a$, we have$$(f\circ g)(a)=\bigvee\limits_{(h,k) 
\in {A_a}} {\min \{ f(h),g(k)\}}\ge \min\{f(u),g(v)\}.$$Since $f$ is 
a fuzzy left ideal of $H$, we have $f(x\circ a)\ge f(a)$ and since 
$u\in x\circ a$, we have $f(u)\ge f(a)$. Since $g$ is a fuzzy left 
ideal of $H$, we have $g(y\circ a)\ge g(a)$ and since $v\in y\circ 
a$, we have $g(v)\ge g(a)$. Hence we obtain$$(f\circ g)(a)\ge 
\min\{f(a),g(a)\}=(f\wedge g)(a).$$$\Longleftarrow$. By Theorem 45, 
it is enough to prove that the left ideals of $H$ are idempotent. Let 
now $A$ be a left ideal of $H$ and $a\in A$. Since $f_A$ is a fuzzy 
left ideal of $H$, by hypothesis, we have $f_A=f_A\wedge f_A\preceq 
f_A\circ f_A$. Then $1=f_A(a)\le f_A\circ f_A)(a)$. If 
$A_a=\emptyset$, then $(f_A\circ f_A)(a)=0$ which is impossible. Thus 
$A_a\not=\emptyset$ and$$(f_A\circ f_A)(a):=\bigvee\limits_{(h,k) \in 
{A_a}} {\min \{ f_A(h),f_A(k)\}}.$$Then there exists $(y,z)\in A_a 
\mbox { such that } y\in A \mbox { and } z\in A\hfill(*)$\\Indeed: 
Suppose there is no $(y,z)\in A_a$ such that $y\in A$ and $z\in A$. 
Then for every $(y,z)\in A_a$ we have $y\notin A$ or $z\notin A$. 
Then for every $(y,z)\in A_a$ we have $f_A(y)=0$ or $f_A(z)=0$. Then 
for every $(y,z)\in A_a$ we have $\min\{f_A(y),f_A(z)\}=0$, so 
$(f_A\circ f_A)(a)=0$ which is impossible. By $(*)$, we have$$a\in 
y\circ z\subseteq A*A\subseteq H*A\subseteq A,$$ thus we have 
$A*A=A$, and $A$ is idempotent. $\hfill\Box$\medskip

\noindent{\bf Theorem 48.} {\it An hypersemigroup $(H,\circ)$ is left 
quasi-regular if and only if the fuzzy left ideals of $H$ are 
idempotent.}\medskip

\noindent{\bf Proof.} $\Longrightarrow$. Let $f$ be a fuzzy left 
ideal of $H$ and $a\in H$. Then $(f\circ f)(a)\le f(a)$. In fact, if 
$A_a=\emptyset$, then $(f\circ f)(a)=0\le f(a)$. Let 
$A_a\not=\emptyset$. Then$$(f\circ f)(a)=\bigvee\limits_{(x,y) \in 
{A_a}} {\min \{ f(x),f(y)\}}.$$On the other hand, 
$$\min\{f(x),f(y)\}\le f(a) \mbox { for every } (x,y)\in 
A_a.$$Indeed: Let $(x,y)\in A_a$. Since $f$ is a fuzzy left ideal of 
$H$ we have $f(x\circ y)\ge f(y)$, and since $a\in x\circ y$, we have 
$f(a)\ge f(y)\ge \min\{f(x),f(y)\}$. Thus we get $(f\circ f)(a)\le 
f(a)$. Moreover $f\preceq f\circ f$. Indeed: Let $a\in H$. Since $H$ 
is left quasi-regular, there exist $x,y\in H$ such that $a\in (x\circ 
a)*(y\circ a)$. Then $a\in u\circ v$ for some $u\in x\circ a$, $v\in 
y\circ a$. Since $(u,v)\in A_a$, we have$$(f\circ 
f)(a):=\bigvee\limits_{(h,k) \in {A_a}} {\min \{ f(h),f(k)\}}\ge 
\min\{f(u),f(v)\}.$$Since $f$ is a fuzzy left ideal of $H$, we have 
$f(x\circ a)\ge f(a)$ and since $u\in x\circ a$, we get $f(u)\ge 
f(a)$. Again since $f$ is a fuzzy left ideal of $H$, we have 
$f(y\circ a)\ge f(a)$ and since $v\in y\circ a$, we get $f(v)\ge 
f(a)$. Thus we have$$(f\circ f)(a)\ge\min\{f(a),f(a)\}=f(a),$$so 
$f\preceq f\circ f$, and $f$ is idempotent.\\$\Longleftarrow$. By 
Theorem 45, it is enough to prove that the left ideals of $H$ are 
idempotent. Let now $A$ be a left ideal of $H$ and $a\in A$. Since 
$f_A$ is a fuzzy left ideal of $H$, by hypothesis, we have 
$f_A=f_A\circ f_A$, thus $1=f_A(a)=(f_A\circ f_A)(a)$. Then $1\le 
(f_A\circ f_A)(a)$ and for the rest of the proof we refer to the 
proof of the $\Leftarrow$-part of the previous theorem. $\hfill\Box$

The concept of right regular semigroups can be naturally transferred 
to hypersemigroups in the following definition.\medskip

\noindent{\bf Definition 49.} An hypersemigroup $(H,\circ)$ is called 
{\it right quasi-regular} if for every $a\in H$ there exist $x,y\in 
H$ such that $a\in (a\circ x)*(a\circ y)$.\medskip

\noindent{\bf Proposition 50.} {\it Let $(H,\circ)$ be an 
hypersemigroup. The following are equivalent:\begin{enumerate}
\item H is right quasi-regular.
\item $a\in \{a\}*H*\{a\}*H$ for every $a\in H$.
\item $A\subseteq A*H*A*H$ for every $A\in {\cal 
P}^*(H)$.\end{enumerate}}

The right analogues of Theorems 43--48 also hold and we have 
the\newline following:\medskip

\noindent{\bf Theorem 51.} {\it Let $(H,\circ)$ be an hypersemigroup. 
The following are equivalent:\begin{enumerate}
\item H is right quasi-regular.
\item $A\cap B\subseteq A*B$ for every nonempty subset $A$ and every 
ideal $B$ of $H$.
\item $A\cap B\subseteq A*B$ for every bi-ideal $A$ and every ideal 
$B$ of $H$.
\item $A\cap B\subseteq A*B$ for every right ideal $A$ and every 
ideal $B$ of $H$.
\item $R(A)\cap I(A)\subseteq R(A)*I(A)$ for every $A\in {\cal 
P}^*(H)$.
\item $R(a)\cap I(a)\subseteq R(a)*I(a)$ for every $a\in H$.
\end{enumerate}}

\noindent{\bf Theorem 52.} {\it An hypersemigroup $(H,\circ)$ is 
right quasi-regular if and only if, for any right ideals $A$ and $B$ 
of $H$, we have $A\cap B\subseteq A*B$.}\medskip

\noindent{\bf Theorem 53.} {\it An hypersemigroup $(H,\circ)$ is 
right quasi-regular if and only if the right ideals of $H$ are 
idempotent.}\medskip

\noindent{\bf Theorem 54.} {\it Let $(H,\circ)$ be an hypersemigroup. 
The following are equivalent:\begin{enumerate}
\item H is right quasi-regular.
\item $f\wedge g\preceq f\circ g$ for every fuzzy subset f and every 
fuzzy ideal g of H.
\item $f\wedge g\preceq f\circ g$ for every fuzzy bi-ideal f and 
every fuzzy ideal g of H.
\item $f\wedge g\preceq f\circ g$ for every fuzzy right ideal f and 
every fuzzy ideal g of H.
  \end{enumerate}}

\noindent{\bf Theorem 55.} {\it An hypersemigroup $(H,\circ)$ is 
right quasi-regular if and only if for any fuzzy right ideals f and g 
of H, we have $f\wedge g\preceq f\circ g$.}

A fuzzy subset $f$ of an hypergroupoid is called {\it idempotent} if 
$f\circ f=f$. \medskip

\noindent{\bf Theorem 56.} {\it An hypersemigroup $(H,\circ)$ is 
right quasi-regular if and only if the fuzzy right ideals of $H$ are 
idempotent.}\medskip

The concept of semisimple semigroups can be naturally transferred to 
hypersemigroups as follows: \medskip

\noindent{\bf Definition 57.} An hypersemigroup $(H,\circ)$ is called 
{\it semisimple} if for every $a\in H$ there exist $x,y,z\in H$ such 
that $a\in (x\circ a)*(y\circ a)*\{z\}$. \medskip

\noindent{\bf Proposition 58.} {\it Let $(H,\circ)$ be an 
hypersemigroup. The following are equivalent:\begin{enumerate}
\item H is semisimple.
\item $a\in H*\{a\}*H*\{a\}*H$ for every $a\in H$.
\item $A\subseteq H*A*H*A*H$ for every $A\in {\cal 
P}^*(H)$.\end{enumerate}}
Let us prove the implication $(3)\Rightarrow (1)$: Let $a\in H$. By 
(3), we have\\$\{a\}\subseteq 
{\bigg(}{\Big(}H*\{a\}{\Big)}*{\Big(}H*\{a\}{\Big)}{\bigg)}*H$. Then 
$$a\in x\circ z \mbox { for some }
x\in {\Big(}H*\{a\}{\Big)}*{\Big(}H*\{a\}{\Big)},\; z\in H.$$Then
$x\in u\circ v$ for some $u,v\in H*\{a\}$, $u\in x\circ a$ for some 
$x\in H$ and $v\in y\circ a$ for some $y\in H$. Then$$a\in x\circ 
z,\, x\in u\circ v,\, u\in x\circ a,\, v\in y\circ a,\, z\in H.$$Thus 
we have\begin{eqnarray*}a\in x\circ z&=&\{x\}*\{z\}\subseteq (u\circ 
v)*\{z\}=\{u\}*\{v\}*\{z\}\\&\subseteq&(x\circ a)*(y\circ 
a)*\{z\}.\end{eqnarray*}Since $x,y,z\in H$ and $a\in (x\circ 
a)*(y\circ a)*\{z\}$, $H$ is semisimple and condition (1) is 
satisfied. $\hfill\Box$\medskip

\noindent{\bf Theorem 59.} {\it Let $(H,\circ)$ be an hypersemigroup. 
The following are equivalent:\begin{enumerate}
\item H is semisimple.
\item The ideals of H are idempotent.
\item $A\cap B=A*B$ for all ideals $A,B$ of H.
\item $I(A)=I(A)*I(A)$ for every $A\in {\cal P}^*(H)$.
\item $I(a)=I(a)*I(a)$ for every $a\in H$.\end{enumerate} }
\noindent{\bf Proof.} $(1)\Longrightarrow (2)$. Let $A$ be an ideal 
of $H$. Since $A\in{\cal P}^*(H)$ and $H$ is semisimple, by 
Proposition 58, we have \begin{eqnarray*}A
&\subseteq&(H*A)*H*(A*H)\subseteq A*H*A=(A*H)*A\subseteq 
A*A\\&\subseteq&A*H\subseteq A,\end{eqnarray*} so $A*A=A$, and $A$ is 
idempotent.\\
$(2)\Longrightarrow (3)$. Let $A, B$ be ideals of $H$. Then 
$A*B\subseteq A*H\subseteq A$ and $A*B\subseteq H*B\subseteq B$, so 
$A*B\subseteq A\cap B$. On the other hand, $A\cap B$ is an ideal of 
$H$ and, by hypothesis, we have $A\cap B=(A\cap B)*(A\cap B)\subseteq 
A*B$. Thus we have $A\cap B=A*B$. \\The implications $(3)\Rightarrow  
(4)$ and $(4)\Rightarrow (5)$ are obvious.\\$(5)\Longrightarrow 
(1)$.
Exactly as in the Lemma 2 in [2], we prove that 
$$I(a)=I(a)*I(a)*I(a)*I(a)*I(a)$$and that 
$$I(a)*I(a)*I(a)*I(a)*I(a)\subseteq H*\{a\}*H*\{a\}*H.$$ Then we get 
$a\in H*\{a\}*H*\{a\}*H$, and $H$ is semisimple. 
$\hfill\Box$\medskip

\noindent{\bf Proposition 60.} {\it Let $(H,\circ)$ be an 
hypersemigroup. Then we have the following:\begin{enumerate}
\item If H is regular, then it is left and right quasi-regular.
\item If H is left (or right) quasi-regular, then it is semisimple.
\item If H is intra-regular, then it is semisimple.\end{enumerate}}
\noindent{\bf Proof.} 1. Let $H$ be regular and $A\in {\cal P}^*(H)$. 
Then we have\begin{eqnarray*}A&\subseteq&A*H*A\subseteq 
A*H*(A*H*A)\subseteq(H*H)*(A*H*A)\\&\subseteq&
H*A*H*A,\end{eqnarray*}so $H$ is left quasi-regular. Similarly $H$ is 
right quasi-regular.\\2. Let $H$ be left quasi-regular and $A\in 
{\cal P}^*(H)$. Then we have\begin{eqnarray*}A&\subseteq& 
H*A*H*A\subseteq 
H*(H*A*H*A)*(H*A)\\&=&(H*H)*(A*H*A)*(H*A)\\&\subseteq&(H*H)*(A*H*A)*(H*H)
\\&\subseteq&H*A*H*A*H,\end{eqnarray*}thus $H$ is semisimple.\\3.
Let $H$ be intra-regular and $A$ a nonempty subset of $H$. Then we 
have\begin{eqnarray*}A&\subseteq &H*A*A*H\subseteq 
H*(H*A*A*H)*A*H\\&=&(H*H)*A*(A*H)*A*H\\&\subseteq&(H*H)*A*(H*H)*A*H\\
&\subseteq&H*A*H*A*H,\end{eqnarray*}and $H$ is semisimple. 
$\hfill\Box$\medskip

\noindent{\bf Theorem 61.} {\it Let H be an hypersemigroup. The 
following are equivalent:\begin{enumerate}
\item H is semisimple.
\item For each fuzzy ideals f and g of H, we have $f\wedge g=f\circ 
g$.
\item For every fuzzy ideal f of H, we have $f=f\circ 
f$.\end{enumerate}}
\noindent{\bf Proof.} $(1)\Longrightarrow (2)$. Let $f$ and $g$ be 
fuzzy ideals of $H$. Since $f$ is a fuzzy right ideal and $g$ is a 
fuzzy left ideal of $H$, we have $f\circ g\preceq f\wedge g$. Let now 
$a\in H$. Then $(f\wedge g)(a)\le (f\circ g)(a)$. In fact: Since $H$ 
is semisimple, there exist $x,y,z\in H$ such that $a\in (x\circ 
a)*(y\circ a)*\{z\}$. Then there exist $u\in x\circ a$ and $v\in 
(y\circ a)*\{z\}$ such that $a\in u\circ v$. Since $v\in (y\circ 
a)*\{z\}$, there exists $w\in y\circ a$ such that $v\in w\circ z$. 
Thus we have$$u\in x\circ a,\; a\in u\circ v,\; w\in y\circ a \mbox { 
and } v\in w\circ z.$$Since $a\in u\circ v$, we have $(u,v)\in A_a$. 
Since $(u,v)\in A_a$, we have$$(f\circ g)(a):=\bigvee\limits_{(h,k) 
\in {A_a}} {\min \{ f(h),g(k)\} }\ge \min\{f(u), g(v)\}.$$Since $f$ 
is a fuzzy left ideal of $H$, we have $f(x\circ a)\ge f(a)$ and since 
$u\in x\circ a$, we have $f(u)\ge f(a)$. Since $g$ is a fuzzy right 
ideal of $H$, we have $g(w\circ z)\ge g(w)$ and since $v\in w\circ 
z$, we have $g(v)\ge g(w)$. Since $g$ is a fuzzy left ideal of $H$, 
we have $g(y\circ a)\ge g(a)$ and since $w\in y\circ a$, we have 
$g(w)\ge g(a)$. Thus we get $g(v)\ge g(a)$. Hence we obtain$$(f\circ 
g)(a)\ge \min\{f(a),g(a)\}=(f\wedge g)(a),$$so $f\wedge g\preceq 
f\circ g$.\\The implication $(2)\Rightarrow (3)$ is 
obvious.\\$(3)\Longrightarrow (1)$. Let $a\in H$. We prove that 
$I(a)\subseteq I(a)*I(a)$. Then, since $I(a)$ is an ideal of $H$, we 
have $I(a)=I(a)*I(a)$ and, by Theorem 59, $H$ is semisimple. Let now 
$b\in I(a)$. Then $b\in I(a)*I(a)$. In fact: Since $I(a)$ is an ideal 
of $H$, the characteristic function $f_{I(a)}$ is a fuzzy ideal of 
$H$. By hypothesis, we have $f_{I(a)}=f_{I(a)}\circ f_{I(a)}$, then 
$f_{I(a)}(b)={\Big(}f_{I(a)}\circ f_{I(a)}{\Big)}(b)$. Since $b\in 
I(a)$, we have $f_{I(a)}(b)=1$, then $1={\Big(}f_{I(a)}\circ 
f_{I(a)}{\Big)}(b)$. If $A_b=\emptyset$, then ${\Big(}f_{I(a)}\circ 
f_{I(a)}{\Big)}(b)=0$ which is no possible. Thus $A_b\not=\emptyset$ 
and
$${\Big(}f_{I(a)}\circ f_{I(a)}{\Big)}(b)=\bigvee\limits_{(y,z) \in 
{A_b}} {\min \{ {f_{I(a)}}(y),{f_{I(a)}}(z)\} }.$$Then there exists 
$(y,z)\in A_b$ such that $y\in I(a)$ and $z\in I(a)$ 
(otherwise,\newline ${\Big(}f_{I(a)}\circ f_{I(a)}{\Big)}(b)=0$ which 
is impossible). Therefore, we have $$b\in y\circ z\subseteq 
I(a)*I(a), \mbox { and  then } b\in I(a)*I(a).$$ $\hfill\Box$

These theorems hold for ordered semigroups as well. Let us give the 
result which corresponds to left quasi-regular ordered 
semigroups.\medskip

\noindent{\bf Definition 62.} An ordered hypersemigroup $(H,\circ)$ 
is called {\it left quasi-regular} if$$\forall\;a\in H \;\,\exists\; 
x,y,t\in H \mbox { such that } t\in (x\circ a)*(y\circ a) \mbox { and 
} a\le t.$$\noindent{\bf Proposition 63.} {\it Let $(H,\circ)$ be an 
hypersemigroup. The following are equivalent:\begin{enumerate}
\item H is left quasi-regular.
\item $a\in {\Big(}H*\{a\}*H*\{a\}{\Big]}$ for every $a\in H$.
\item $A\subseteq (H*A*H*A]$ for every $A\in {\cal 
P}^*(H)$.\end{enumerate}}
\noindent{\bf Theorem 64.} {\it Let $(H,\circ,\le)$ be an ordered 
hypersemigroup. The following are equivalent:\begin{enumerate}\item 
$H$ is left quasi-regular.\item $A\cap B\subseteq (A*B]$ for every 
ideal A and every nonempty subset B of H.\item $A\cap B\subseteq 
(A*B]$ for every ideal A and every bi-ideal B of H.\item $A\cap 
B\subseteq (A*B]$ for every ideal A and eery left ideal B of H.\item 
$I(A)\cap L(A)\subseteq {\Big(}I(A)\cap L(A){\Big]}$ for every $A\in 
{\cal P}^*(H)$.\item $I(a)\cap L(a)\subseteq {\Big(}I(a)\cap 
L(a){\Big]}$ for every $a\in H$.\end{enumerate}}
\noindent{\bf Proof.} $(1)\Longrightarrow (2)$. Let $A$ be an ideal, 
$B$ a nonempty subset of $H$ and $a\in A\cap B$. Since $a\in H$ and 
$H$ is left quasi-regular, we 
have$$a\in{\Big(}H*\{a\}*H*\{a\}{\Big]}\subseteq
{\Big(}(H*A*H)*B{\Big]}\subseteq (A*B].$$The implications 
$(2)\Rightarrow (3)$ and $(2)\Rightarrow (3)$ and $(4)\Rightarrow 
(5)\Rightarrow (6)$ are obvious and $(3)\Rightarrow (4)$ since every 
left ideal of $H$ is a bi-ideal of $H$.\\$(6)\Longrightarrow (1)$. 
Let $a\in H$. By hypothesis, we have\begin{eqnarray*}a&\in&
{\Big(}I(a)*L(a){\Big]}\\&=&{\Bigg(}{\bigg(}\{a\}\cup 
{\Big(}H*\{a\}{\Big)}\cup {\Big(}\{a\}*H{\Big)}\cup 
{\Big(}H*\{a\}*H{\Big)}{\bigg]}*{\bigg(}\{a\}\cup 
{\Big(}H*\{a\}{\Big)}{\bigg]}{\Bigg]}\\&=&{\Bigg(}{\bigg(}\{a\}\cup 
{\Big(}H*\{a\}{\Big)}\cup {\Big(}\{a\}*H{\Big)}\cup 
{\Big(}H*\{a\}*H{\Big)}{\bigg)}*{\bigg(}\{a\}\cup 
{\Big(}H*\{a\}{\Big)}{\bigg)}{\Bigg]}\\&=&{\bigg(}{\Big(}\{a\}*\{a\}{\Big)}\cup 
{\Big(}H*\{a\}*\{a\}{\Big)}\cup {\Big(}\{a\}*H*\{a\}{\Big)}\cup 
{\Big(}H*\{a\}*H*\{a\}{\Big)}{\bigg]}\\&=&{\Big(}\{a\}*\{a\}{\Big]}\cup 
{\Big(}H*\{a\}*\{a\}{\Big]}\cup {\Big(}\{a\}*H*\{a\}{\Big]}\cup 
{\Big(}H*\{a\}*H*\{a\}{\Big]}. \end{eqnarray*}If $a\in 
{\Big(}\{a\}*\{a\}{\Big]}$, then we 
have\begin{eqnarray*}\{a\}&\subseteq&{\Big(}\{a\}*\{a\}{\Big]}\subseteq 
{\bigg(}{\Big(}\{a\}*\{a\}{\Big]}*{\Big(}\{a\}*\{a\}{\Big]}{\bigg]}
\\&=&{\bigg(}{\Big(}\{a\}*\{a\}{\Big)}*{\Big(}\{a\}*\{a\}{\Big)}{\bigg]}
\subseteq {\Big(}H*\{a\}*H*\{a\}{\Big]},\end{eqnarray*}so $a\in 
{\Big(}H*\{a\}*H*\{a\}{\Big]}$. If $a\in 
{\Big(}H*\{a\}*\{a\}{\Big]}$, then\begin{eqnarray*}a\in 
\{a\}&\subseteq&{\bigg(}H*{\Big(}H*\{a\}*\{a\}{\Big]}*\{a\}{\bigg]}=
{\bigg(}(H]*{\Big(}H*\{a\}*\{a\}{\Big]}*{\Big(}\{a\}{\Big]}{\bigg]}\\&=&
{\bigg(}H*{\Big(}H*\{a\}*\{a\}{\Big)}*\{a\}{\bigg]}\subseteq 
{\Big(}H*\{a\}*H*\{a\}{\Big]}.\end{eqnarray*}If $a\in 
{\Big(}\{a\}*H*\{a\}{\Big]}$, in a similar way we get $a\in 
{\Big(}H*\{a\}*H*\{a\}{\Big]}$. In each case $H$ is left 
quasi-regular. $\hfill\Box$ \medskip

\noindent{\bf Theorem 65.} {\it An ordered hypersemigroup 
$(H,\circ,\le)$ is left quasi-regular if and only if, for any left 
ideals A, B of H, we have $A\cap B\subseteq (A*B]$.}\smallskip

\noindent{\bf Proof.} $\Longrightarrow$. Let $A$, $B$ be left ideals 
of $H$ and $a\in A\cap B$. Since $a\in H$ and $H$ is left 
quasi-regular, there exist $x,y,t\in H$ such that $t\in (x\circ 
a)*(y\circ a)$ and $a\le t$. We have $a\le t\in (x\circ a)*(y\circ 
a)\in (H*A)*(H*B)\subseteq A*B$, so $a\in (A*B]$.\\$\Longleftarrow$. 
Let $A\in{\cal P}^*(H)$. Then $A\subseteq (H*A*H*A]$. Indeed: By 
hypothesis, we have\begin{eqnarray*}A&\subseteq&L(A)\cap 
L(A)\subseteq (L(A)*L(A)]={\bigg(}{\Big(}A\cup 
(H*A]{\Big)}*{\Big(}A\cup 
(H*A]{\Big)}{\bigg)}\\&\subseteq&{\bigg(}{\Big(}(A]\cup 
(H*A]{\Big)}*{\Big(}(A]\cup 
(H*A]{\Big)}{\bigg)}\\&=&{\bigg(}{\Big(}(A]*(A]{\Big)}\cup 
{\Big(}(H*A]*(A]{\Big)}\cup {\Big(}(A]*(H*A]{\Big)}\cup 
{\Big(}(H*A]*(H*A]{\Big)}{\bigg]}\\&\subseteq&{\Big(}(A*A]\cup 
(H*A*A]\cup (A*H*A]\cup (H*A*H*A]{\Big]}. 
\end{eqnarray*}Then\begin{eqnarray*}A*A&\subseteq&(A]*{\Big(}(A*A]\cup 
(H*A*A]\cup (A*H*A]\cup 
(H*A*H*A]{\Big]}\\&=&{\Big(}(A]*(A*A]{\Big)}\cup 
{\Big(}(A]*(H*A*A]{\Big)}\cup {\Big(}(A]*(A*H*A]{\Big)}\\
&\cup& {\Big(}(A]*(H*A*H*A]{\Big)}\\&\subseteq&(A*A*A]\cup 
(A*H*A*A]\cup (A*A*H*A]\cup (A*H*A*H*A]\\&\subseteq& 
(A*H*A],\end{eqnarray*}and so $(A*A]\subseteq 
{\Big(}(A*H*A]{\Big]}=(A*H*A]$. Then$$H*A*A\subseteq 
(H]*(A*H*A]\subseteq (H*A*H*A],$$from which $(H*A*A]\subseteq 
(H*A*H*A]$. Thus we get
$$A\subseteq {\Big(}(A*H*A]\cup 
(H*A*H*A]{\Big]}.$$Then\begin{eqnarray*} 
A*H*A&\subseteq&{\Big(}(A*H*A]\cup 
(H*A*H*A]{\Big]}*H*A\\&\subseteq&{\Big(}(A*H*A]\cup 
(H*A*H*A]{\Big]}*(H]*(A]\\&\subseteq&{\Big(}(A*H*A]\cup 
(H*A*H*A]{\Big]}*(H*A]\\&\subseteq&{\bigg(}{\Big(}(A*H*A]\cup 
(H*A*H*A]{\Big)}*H*A{\bigg]}\\&=&{\bigg(}{\Big(}(A*H*A]*H*A{\Big)}\cup 
{\Big(}(H*A*H*A]*H*A{\Big)}{\bigg]}\\&\subseteq&{\bigg(}{\Big(}(A*H*A]*(H*A]{\Big)}\cup 
{\Big(}(H*A*H*A]*(H*A]{\Big)}{\bigg]}\\&\subseteq&{\Big(}(A*H*A*H*A]\cup 
(H*A*H*A*H*A]{\Big]}\\&\subseteq&{\Big(}(H*A*H*A]{\Big]}=(H*A*H*A].\end{eqnarray*}Then 
$A\subseteq (H*A*H*A]$, and $H$ is left quasi-regular. 
$\hfill\Box$\medskip

\noindent{\bf Theorem 66.} {\it An ordered hypersemigroup 
$(H,\circ,\le)$ is left quasi-regular if and only if the left ideals 
of H are idempotent}.\medskip

\noindent{\bf Proof.} $\Longrightarrow.$ Let $L$ be a left ideal of 
$H$. By Theorem 65, we have $L=(L*L]\subseteq (H*L]\subseteq (L]=L$, 
so $(L*L]=L$.\\$\Longleftarrow.$ By Theorem 64(4)$\Rightarrow (1)$, 
it is enough to prove that for every ideal $A$ and every left ideal 
$B$ of $H$, we have $A\cap B\subseteq (A*B]$. Let now $A$ be an ideal 
and $B$ a left ideal of $H$. Then $A\cap B$ is an ideal of $H$. 
Indeed: Take an element $a\in A$ and an element $b\in B$ 
$(A,B\not=\emptyset)$. Then $a\circ b\subseteq A*B\subseteq 
A*H\subseteq A$ and $a\circ b\subseteq A*B\subseteq H*B\subseteq B$, 
so $a\circ b\subseteq A*B$. Since $a\circ b\not=\emptyset$, the set 
$A*B$ is a nonempty subset of $H$. In addition, we have $H*(A\cap 
B)\subseteq (H*A)\cap (H*B)\subseteq A\cap B$ and if $x\in A\cap B$ 
and $H\ni y\le x$ then, since $y\le x\in A$, we have $y\in A$ and 
since $y\le x\in B$ we have $y\in B$, thus we have $y\in A\cap B$. 
Since $A\cap B$ is a left ideal of $H$, by hypothesis, we have $A\cap 
B={\Big(}(A\cap B)*(A\cap B){\Big]}\subseteq (A*B]$. $\hfill\Box$

For an ordered hypergroupoid $H$, we denote by $A_a$ the relation on 
$H$ defined by $A_a:=\{(y,z) \mid \exists \; u\in y\circ z \mbox { 
such that } a\le u\}$.

For an ordered hypersemigroup, let us prove the theorem which 
corresponds to the Theorem 47. \medskip

\noindent{\bf Theorem 67.} {\it An ordered hypersemigroup 
$(H,\circ,\le)$ is left quasi-regular if and only if, for any fuzzy 
left ideals f, g of H, we have $f\wedge g\preceq f\circ g$.}\medskip

\noindent{\bf Proof.} $\Longrightarrow$. Let $f$ and $g$ be fuzzy 
left ideals of $H$ and $a\in H$. Then $(f\wedge g)(a)\le (f\circ 
g)(a)$. Indeed: Since $a\in H$ and $H$ is left quasi-regular, there 
exist $x,y,t\in H$ such that $t\in (x\circ a)*(y\circ a)$ and $a\le 
t$. Since $a\le t$ and $t\in u\circ v$, we have $(u,v)\in A_a$. Then 
we have$$(f\circ g)(a)=\bigvee\limits_{(h,k) \in {A_a}} {\min \{ 
f(h),g(k)\} }\ge \min\{f(u),g(v)\}.$$Since $f$ is a fuzzy left ideal 
of $H$, we have $f(x\circ a)\ge f(a)$ and, since $u\in x\circ a$, we 
have $f(u)\ge f(a)$. Since $g$ is a fuzzy left ideal of $H$, we have 
$g(y\circ a)\ge g(a)$ and, since $v\in y\circ a$, we have $g(v)\ge 
g(a)$. Thus we have $(f\circ g)(a)\ge \min\{f(a),g(a)\}=(f\wedge 
g)(a)$.\\$\Longleftarrow$. By Theorem 66, it is enough to prove that 
the left ideals of $H$ are idempotent. Let now $A$ be a left ideal of 
$H$ and $a\in A$. Since $f_A$ is a fuzzy left ideal of $H$, by 
hypothesis, we have $f_A=f_A\wedge f_A\preceq f_A\circ f_A$. Then 
$1=f_A(a)\le (f_A\circ f_A)(a)$. If $A_a=\emptyset$, then $(f_A\circ 
f_A)(a)=0$ which is impossible. Thus $A_a\not=\emptyset$ 
and$$(f_A\circ f_A)(a):=\bigvee\limits_{(h,k) \in {A_a}} {\min \{ 
f_A(h),f_A(k)\}}.$$Then there exists $(y,z)\in A_a \mbox { such that 
} y\in A \mbox { and } z\in A$.\\Since $(y,z)\in A_a$, there exists 
$u\in y\circ z$ such that $a\le u$. Since $y,z\in A$, we have $y\circ 
z\subseteq A*A$, so $u\in A*A$. Since $a\le u\in A*A$, we have $a\in 
(A*A]$. We also have $(A*A]\subseteq (H*A]\subseteq (A]=A$, hence we 
obtain $A=(A*A]$ and the proof is completed. $\hfill\Box$

In the next theorem we give the analogous of the Theorems 46 and 48 
in case of ordered hypersemigroup.\medskip

\noindent{\bf Theorem 68.} {\it Let $(H,\circ,\le)$ be an ordered 
hypersemigroup. The following are equivalent:\begin{enumerate}
\item H is left quasi-regular.
\item $f\wedge g\preceq f\circ g$ for every fuzzy ideal f and every 
fuzzy subset g of H.
\item $f\wedge g\preceq f\circ g$ for every fuzzy ideal f and every 
fuzzy bi-ideal g of H.
\item $f\wedge g\preceq f\circ g$ for every fuzzy ideal f and every 
fuzzy left ideal g of H.
\item The fuzzy left ideals of $H$ are idempotent.\end{enumerate}}

\noindent{\bf 4.} If $S$ is an semigroup, then the relation ``$\cal 
N$" on $S$ defined by $x{\cal N}y$ $\Leftrightarrow$ $N(x)=N(y)$,  
where $N(a)$ is the filter of $S$ generated by $a$ $(a\in S)$, is the 
least semilattice congruence on $S$. If $S$ is an ordered semigroup, 
then ${\cal N}$ is not the least semilattice congruence on $S$ in 
general, but it is the least complete semilattice congruence on $H$ 
[19]. We examine the same for hypersemigroups.

The concepts of congruences and semilattices congruences of 
semigroups are naturally transferred to hypersemigroups as 
follows:\medskip

\noindent{\bf Definition 69.} Let $H$ be an hypergroupoid. An 
equivalence relation $\sigma$ on $H$ is called {\it right congruence} 
if $(a,b)\in\sigma$ implies $(a\circ c, b\circ c)\in\sigma$ for every 
$c\in H$, in the sense that if $u\in a\circ c$ and $v\in b\circ c$, 
then $(v,v)\in\sigma$. It is called {\it left congruence} if 
$(a,b)\in\sigma$ implies $(c\circ a, c\circ b)\in\sigma$ for every 
$c\in H$, in the sense that if $u\in c\circ a$ and $v\in c\circ b$, 
then $(u,v)\in\sigma$. By a {\it congruence} on $H$ we mean a 
relation on $H$ which is both a right and a left congruence on 
$H$.\medskip

\noindent{\bf Definition 70.} Let $H$ be an hypergroupoid. A 
congruence $\sigma$ on $H$ is called {\it semilattice congruence} if, 
for any $a,b\in H$, we have$$(a\circ a, a)\in\sigma \mbox { and } 
(a\circ b, b\circ a)\in\sigma$$that is, for every $a\in H$ and every 
$u\in a\circ a$, we have $(u,a)\in\sigma$ and for every $a,b\in H$, 
and every $u\in a\circ b$ and $v\in b\circ a$, we have 
$(u,v)\in\sigma$.

The concept of a filter of groupoids is naturally transferred to 
hypergroupoids as follows: \medskip

\noindent{\bf Definition 71.} Let $H$ be an hypergroupoid. A nonempty 
subset $F$ of $H$ is called a {\it filter} of $H$ if the following 
assertions are satisfied:\begin{enumerate}
\item if $x,y\in F$, then $x\circ y\subseteq F$
\item if $x,y\in H$ such that $x\circ y\subseteq F$, then $x\in F$ 
and $y\in F$
\item for any $x,y\in H$, we have $x\circ y\subseteq F$ or $(x\circ 
y)\cap F=\emptyset$.\end{enumerate}

For an element $x$ of $H$, we denote by $N(x)$ the filter of $H$ 
generated by $x$, and by $\cal N$ the equivalence relation on $H$ 
defined by$${\cal N}:=\{(x,y)\in H\times H \mid N(x)=N(y)\}.$$

As a modification of the property 1) in [10, the Proposition], we 
have the following proposition\medskip

\noindent{\bf Proposition 72.} {\it If $H$ is an hypergroupoid, then 
the equivalence relation $\cal N$ is a semilattice congruence on 
$H$.}\medskip

\noindent{\bf Proof.} Let $(x,y)\in\cal N$ and $z\in H$. Then 
($(z\circ x, z\circ y)\in {\cal N}$. In fact: Let $u\in z\circ x$ and 
$v\in z\circ y$. Then $(u,v)\in \cal N$. Indeed: Since $u\in N(u)$ 
and $u\in z\circ x$, we have $(z\circ x)\cap N(u)\not=\emptyset$. 
Since $N(u)$ is a filter of $H$, we have $z\circ x\subseteq N(u)$, 
then $z,x\in N(u)$. Since $x\in N(u)$, we have $N(x)\subseteq N(u)$, 
then $y\in N(u)$. Since $z,y\in N(u)$, we have $z\circ y\subseteq 
N(u)$, so $v\in N(u)$, and $N(v)\subseteq N(u)$. By symmetry, we get 
$N(u)\subseteq N(v)$, so we have $N(u)=N(v)$, and $(u,v)\in\cal N$. 
Thus $\cal N$ is a left congruence on $H$. In a similar way we prove 
that $\cal N$ is a right congruence on $H$, and so $\cal N$ is a 
congruence on $H$. Let now $x\in H$. Then $(x\circ x, x)\in {\cal 
N}$. In fact: Let $u\in x\circ x$. Then $(u,x)\in \cal N$. Indeed: 
Since $u\in N(u)$, we have $(x\circ x)\cap N(u)\not=\emptyset$. Since 
$N(u)$ is a filter of $H$, we have $x\circ x\subseteq N(u)$, then 
$x\in N(u)$, and $N(x)\subseteq N(u)$. On the other hand, since $x\in 
N(x)$ and $N(x)$ is a filter of $H$, we have $x\circ x\subseteq 
N(x)$. Then $u\in N(x)$, so $N(u)\subseteq N(x)$. Thus we have 
$N(u)=N(x)$, and $(u,x)\in\cal N$. Let $x,y\in H$. Then $(x\circ y, 
y\circ x)\in {\cal N}.$ In fact: Let $u\in x\circ y$ and $v\in y\circ 
x$. Then $(u,v)\in\cal N$. Indeed, since $u\in N(u)$, we have 
$(x\circ y)\cap N(u)\not=\emptyset$, then $x\circ y\subseteq N(u)$, 
$x,y\in N(u)$, and $y\circ x\subseteq N(u)$. Then $v\in N(u)$, and 
$N(v)\subseteq N(u)$. By symmetry, we get $N(u)\subseteq N(v)$, then 
$N(u)=N(v)$, and $(u,v)\in \cal N$. $\hfill\Box$

The concepts of prime subsets of groupoids can be naturally 
transferred to hypersemigroups as follows. It should be mentioned 
here that Kehayopulu uses the terms ``prime", ``weakly prime" instead 
of ``completely prime", ``prime" considered by Petrich.\medskip

\noindent{\bf Definition 73.} Let $H$ be an hypergroupoid. A nonempty 
subset $T$ of $H$ is called a {\it prime subset} of $H$ if the 
following assertions are satisfied:\begin{enumerate}
\item if $a,b\in H \mbox { such that } a\circ b\subseteq T \mbox { 
then } a\in T \mbox { or } b\in T$ and
\item for every $a,b\in H$, we have $a\circ b\subseteq T$ or $(a\circ 
b)\cap T=\emptyset$.\end{enumerate}
\noindent{\bf Notation 74.} For a subset $I$ of $H$, we denote by 
$\sigma_I$ the equivalence relation on $H$ defined 
by:$$\sigma_I:=\{(a,b)\in H\times H \mid a, b\in I \mbox { or } a, 
b\notin I\}$$(i.e. $a, b$ both belong to $I$ or $a, b$ both do not 
belong t0 $I$).\medskip

\noindent{\bf Lemma 75.} {\it Let H be an hypergroupoi, I a subset of 
H and $a,b,c\in H$. Then

1. if $a\circ c, b\circ c\subseteq I$, then $(a\circ c, b\circ 
c)\in\sigma_I$.\\
Suppose now that, for every $a,b\in H$, we have $a\circ b\subseteq I$ 
or $(a\circ b)\cap I=\emptyset$. Then

2. If $a\circ c, b\circ c\nsubseteq I$, then $(a\circ c, b\circ 
c)\in\sigma_I$.}\medskip

\noindent{\bf Proof.} (1) Let $a\circ c, b\circ c\subseteq I$ and let 
$u\in a\circ c$, $v\in b\circ c$. Then $u,v\in I$, so 
$(u,v)\in\sigma_I$. (2) Let $a\circ c, b\circ c\nsubseteq I$, $u\in 
a\circ c$ and $v\in b\circ c$. If $u,v\in I$, then 
$(u,v)\in\sigma_I$. Let $u\notin I$. Then $v\notin I$. Indeed: 
Suppose $v\in I$. Since $(b\circ c)\cap I\not=\emptyset$, by 
hypothesis, we have $b\circ c\subseteq I$ which is impossible. So 
$u,v\notin I$, and $(u,v)\in\sigma_I$. If $v\notin I$, in a similar 
way we get $u\notin I$, so again $(u,v)\in\sigma_I$. 
$\hfill\Box$\medskip

\noindent{\bf Proposition 76.} (see the analogous for ordered 
semigroups in [10]) {\it Let H be an hypergroupoid and I a prime 
ideal of H. Then the equivalent relation $$\sigma_I:=\{(a,b)\in 
H\times H \mid a,b\in I \mbox { or } a,b\not\in I\}$$ is a 
semilattice congruence on H.}\medskip

\noindent{\bf Proof.} (1) $\sigma_I$ is a congruence on $H$. In fact: 
Let $(a,b)\in\sigma_I$ and $c\in H$. Then $(a\circ c, b\circ 
c)\in\sigma_I.$ In fact: Since $(a,b)\in\sigma_I$, we have $a,b\in I$ 
or $a,b\notin I$. Let $a,b\in I$. Since $I$ is an ideal of $H$, we 
have $a\circ c, b\circ c\subseteq I$. Then, by Lemma 75(1), we have 
$(a\circ c, b\circ c)\in\sigma_I$. Let $a,b\notin I$. Since $c\in H$, 
we have $c\in I$ or $c\notin I$. Let $c\in I$. Since $I$ is an ideal 
of $H$, we have $a\circ c, b\circ c\subseteq I$, then $(a\circ c, 
b\circ c)\in\sigma_I$. Let $c\notin I$. Then $a\circ c, b\circ 
c\nsubseteq I$. Indeed: If $a\circ c\subseteq I$ then, since $I$ is a 
prime ideal of $H$, we have $a\in I$ or $c\in I$, then $a\in I$ which 
is impossible. If $b\circ c\subseteq I$, then $b\in I$ or $c\in I$, 
then $b\in I$ which again is not possible. Since $a\circ c, b\circ 
c\nsubseteq I$, by Lemma 75(2), we have $(a\circ c, b\circ 
c)\in\sigma_I$. Thus $\sigma_I$ is a right congruence on $M$. In a 
similar way we can prove that $\sigma_I$ is a left congruence on $M$, 
so it is a congruence on $M$.

(2) The congruence $\sigma_I$ is a semilattice congruence on $H$. In 
fact:\\Let $a\in H$. Then $(a\circ a,a)\in\sigma_I$. Indeed: Let 
$u\in a\circ a$. Clearly, $a\in I$ or $a\notin I$. If $a\in I$ then, 
since $I$ is an ideal of $H$, we have $a\circ a\subseteq I$, then 
$u\in I$. Since $u,a\in I$, we have $(u,a)\in\sigma_I$. Let $a\notin 
I$. Then $a\circ a\nsubseteq I$. Indeed: If $a\circ a\subseteq I$ 
then, since $I$ is prime, we have $a\in I$ which is impossible. Since 
$a\notin I$ and $a\circ a\nsubseteq I$, we have $\{a\}, a\circ 
a\nsubseteq I$. Then we have $(\{a\}, a\circ , a)\in\sigma_I$, which 
means that $(a, a\circ a)\in \sigma_I$. Let $a,b\in H$. Then  
$(a\circ b, b\circ a)\in\sigma_I$. In fact: If $a\circ b\subseteq I$ 
then, since $I$ is prime, we have $a\in I$ or $b\in I$. Since $I$ is 
an ideal of $H$, we have $b\circ a\subseteq I$. Since $a\circ b, 
b\circ a\subseteq I$, by Lemma 75(1), we have $(a\circ b, b\circ 
a)\in\sigma_I$. If $a\circ b\nsubseteq I$, then $b\circ a\nsubseteq 
I$. Indeed: If $b\circ a\subseteq I$ then, since $I$ is prime, we 
have $b\in I$ or $a\in I$ and, since $I$ is an ideal of $H$, we have 
$a\circ b\subseteq I$ which is impossible. Since $a\circ b, b\circ 
a\nsubseteq I$, by Lemma 75(2), we have $(a\circ b, b\circ 
a)\in\sigma_I$. $\hfill\Box$\medskip

\noindent{\bf Lemma 77.} {\it Let H be an hypergroupoid, A a nonempty 
subset of H, $x, z\in H$. Then we have the 
following:\begin{enumerate}
\item if $\sigma$ is a right congruence on H and $(x, A)\in\sigma$, 
then ${\Big(}x\circ z, A*\{z\}{\Big)}\in\sigma$.
\item if $\sigma$ a left congruence on H and $(x, A)\in\sigma$, then 
${\Big(}z\circ x, \{z\}*A{\Big)}\in\sigma$.
\item if $\sigma$ is a right congruence on H and $(A,B)\in\sigma$, 
then ${\Big(}A*\{x\},B*\{x\}{\Big)}\in\sigma$.
\item if $\sigma$ is a left congruence on H and $(A,B)\in\sigma$, 
then ${\Big(}\{x\}*A,\{x\}*B{\Big)}\in\sigma$.\end{enumerate}}
\noindent{\bf Proof.} (1) Let $\sigma$ be a right congruence on $H$, 
$(x, A)\in\sigma$, $u\in x\circ z$ and $v\in A*\{z\}$. Then 
$(u,v)\in\sigma$. Indeed: Since $v\in A*\{z\}$, we have $v\in a\circ 
z$ for some $a\in A$. Since $(x, A)\in\sigma$ and $a\in A$, we have 
$(x,a)\in\sigma$. Since $\sigma$ is a right congruence on $H$, we 
have $(x\circ z, a\circ z)\in\sigma$. Since $u\in x\circ z$ and $v\in 
a\circ z$, we have $(u,v)\in\sigma$.

(3) Let $\sigma$ be a right congruence on $H$ and $(A,B)\in\sigma$.
Let $u\in A*\{x\}$, $v\in B*\{x\}$. Then $u\in a\circ x$ for some 
$a\in A$ and $v\in b\circ x$ for some $b\in B$. Since 
$(A,B)\in\sigma$, $a\in A$ and $b\in B$, we have $(a,b)\in\sigma$. 
Since $\sigma$ is a right congruence on $H$, we have $(a\circ 
x,b\circ x)\in\sigma$. Since $u\in a\circ x$, $v\in b\circ x$, we 
have $(u,v)\in\sigma$.\\The properties (2) and (4) can be proved in a 
similar way. $\hfill\Box$\medskip

\noindent{\bf Lemma 78.} {\it If H is an hypergroupoid B a nonempty 
subset of H, $(x,A)\in\sigma$ and $B\subseteq A$, then 
$(x,B)\in\sigma$.}\smallskip

\noindent{\bf Proof.} Let $b\in B$. Since $b\in A$ and $(x,A)\in 
\sigma$, we have $(x,b)\in\sigma$.$\hfill\Box$

As a modification of the Lemma in [10], we have the following 
proposition\medskip

\noindent{\bf Proposition 79.} {\it Let H be an hypergroupoid and 
$\emptyset\not=F\subseteq H$. The following are equivalent:

1. F is a filter of H.

2. either $H\backslash F=\emptyset$ or $H\backslash F$ is a prime 
ideal of $H$.}\medskip

\noindent{\bf Proof.} $(1)\Longrightarrow (2)$. Let $H\backslash 
F\not=\emptyset$. Then the set $H\backslash F$ is an ideal of $H$. In 
fact: Let $x\in H$ and $y\in H\backslash F$. Then $x\circ y\subseteq 
H\backslash F$. Indeed: Let $x\circ y\nsubseteq H\backslash F$. Then 
there exists $z\in x\circ y$ such that $z\notin H\backslash F$. Since 
$z\in F$ and $z\in x\circ y$, we have $z\in (x\circ y)\cap F$. Since 
$F$ is a filter of $H$ and $(x\circ y)\cap F\not=\emptyset$,  we have 
$x\circ y\subseteq F$ which is impossible. Thus $S\backslash F$ is a 
left ideal of $H$. Similarly, $H\backslash F$ is a right ideal of 
$H$. The ideal $H\backslash F$ is prime. In fact: (a) Let $x,y\in H$ 
such that $x\circ y\subseteq H\backslash F$. Then $x\in H\backslash 
F$ or $y\in H\backslash F$. Indeed, if $x\in F$ and $y\in F$ then, 
since $F$ is a subgroupoid of $H$, we have $x\circ y\subseteq F$ 
which is impossible. (b) Let $x,y\in H\backslash F$ such that 
$(x\circ y)\cap H\backslash F\not=\emptyset$. Then $x\circ y\subseteq 
H\backslash F$. Indeed: Let $a\in x\circ y$. Suppose $a\not\in 
H\backslash F$. Then $a\in F$. Since $a\in (x\circ  y)\cap F$ and $F$ 
is a filter of $H$, we have $x\circ y\subseteq F$. Then $(x\circ 
y)\cap H\backslash F\subseteq F\cap H\backslash F=\emptyset$, so 
$(x\circ y)\cap H\backslash F=\emptyset$ which is 
impossible.\smallskip

\noindent$(2)\Longrightarrow (1)$. If $H\backslash F=\emptyset$, then 
$F=H$, so $F$ is a filter of $H$. Let $H\backslash F$ be a prime 
ideal of $H$. Then we have the following: (a) Let $x,y\in F$. Then 
$x\circ y\subseteq F$. Indeed: Let $z\in x\circ y$. If $z\notin F$, 
then $z\in H\backslash F$, so $z\in (x\circ y)\cap H\backslash F$. 
Since $(x\circ y)\cap H\backslash F\not=\emptyset$ and $H\backslash 
F$ is prime, we have $x\circ y\subseteq H\backslash F$. Again since 
$H\backslash F$ is prime, we have $x\in H\backslash F$ or $y\in 
H\backslash F$ which is impossible. Thus we have $z\in F$, and $F$ is 
a subgroupoid of $H$. (b) Let $x,y\in H$ such that $x\circ y\subseteq 
F$. Then $x\in F$ and $y\in F$. In fact: If $x\in H\backslash F$, 
then $x\circ y\subseteq (H\backslash F)*H\subseteq H\backslash F$ 
which is impossible. If $y\in S\backslash F$, we also get a 
contradiction. So we have $x,y\in F$. (c) Let $x,y\in H$ such that 
$(x\circ y)\cap F\not=\emptyset$. Then $x\circ y\subseteq F$. Indeed: 
Let $x\circ y\nsubseteq F$. Then there exists $z\in x\circ y$ such 
that $z\not\in F$, so $z\in (x\circ y)\cap H\backslash F$. Since 
$(x\circ y)\cap H\backslash F\not=\emptyset$ and $H\backslash F$ is 
prime, we have $x\circ y\subseteq H\backslash F$. Then $(x\circ 
y)\cap F\subseteq H\backslash F\cap F=\emptyset$, so $(x\circ y)\cap 
F=\emptyset$ which is impossible. $\hfill\Box$\medskip

\noindent{\bf Lemma 80.} {\it Let $(H,\circ)$ be an hypergroupoid, 
${\cal P}^*(H)$. Then\begin{enumerate}
\item if $\sigma$ is a symmetric relation on $(H,\circ)$, then it is 
symmetric on ${\cal P}^*(H)$.
\item if $\sigma$ is a transitive relation on $(H,\circ)$, then it is 
transitive on ${\cal P}^*(H)$.\end{enumerate}}
\noindent{\bf Proof.} 1. Let $(A,B)\in\sigma$, $b\in B$ and $a\in A$. 
Then $(a,b)\in\sigma$ and, since $\sigma$ is symmetric, we have 
$(b,a)\in\sigma$. Thus we have $(B,A)\in\sigma$.

2. Let $(A,B)\in\sigma$, $(B,C)\in\sigma$, $a\in A$ and $c\in C$. 
Take an element $b\in B$ $(B\not=\emptyset)$. We have 
$(a,b)\in\sigma$ and $(b,c)\in\sigma$ and, since $\sigma$ is a 
transitive relation on $H$, we have $(a,c)\in\sigma$. Thus we have 
$(A,C)\in\sigma$. $\hfill\Box$\\An ideal $I$ of $H$ is called {\it 
proper} if $I\not=H$.\medskip

\noindent{\bf Proposition 81.} {\it Let H be an hypersemigroup and 
$\sigma$ be a semilattice congruence on H. Then there exists a family 
$\cal A$ of proper prime ideals of H such that$$\sigma= 
\bigcap\limits_{I \in \cal A} {{\sigma _{\rm I}}}.$$}{\bf Proof.} Let 
$x\in H$. We consider the set$$A_x:=\{y\in H \mid (x, x\circ 
y)\in\sigma\}.$$The set $A_x$ is a filter of $H$. In fact: Since 
$x\in H$ and $\sigma$ is a semilattice congruence on $H$, we have 
$(x, x\circ x)\in\sigma$, so $x\in A_x$, that is, $A_x$ is a nonempty 
subset of $H$. Let $y,z\in A_x$. Then $y\circ z\subseteq A_x$. In 
fact: Let $u\in y\circ z$. Then $u\in A_x$, that is $(x, x\circ 
u)\in\sigma$. Indeed: Let $v\in x\circ u$. By hypothesis, we have 
$(x,x\circ y)\in\sigma$ and $(x,x\circ z)\in\sigma$. Since $(x,x\circ 
y)\in\sigma$, by Lemma 77(1), we have ${\Big(}x\circ z, (x\circ 
y)*\{z\}{\Big)}\in\sigma$. Since $(x, x\circ z)\in\sigma$ and 
${\Big(}x\circ z, (x\circ y)*\{z\}{\Big)}\in\sigma$, by Lemma 80(2), 
we have
${\Big(}x, (x\circ y)*\{z\}{\Big)}\in\sigma$, that is ${\Big(}x, 
\{x\}*(y\circ z){\Big)}\in\sigma$. Since $u\in y\circ z$, we have 
$\{u\}\subseteq y\circ z$, then$$v\in x\circ u=\{x\}*\{u\}\subseteq 
\{x\}*(y\circ z).$$Since ${\Big(}x, \{x\}*(y\circ z){\Big)}\in\sigma$ 
and $v\in\{x\}*(y\circ z)$, we have $(x,v)\in\sigma$.

Let $y,z\in H$ and $y\circ z\subseteq A_x$. Then $y\in A_x$ and $z\in 
A_x$. In fact: \\Since $y\circ z\subseteq A_x$, we have ${\Big(}x, 
\{x\}*(y\circ z){\Big)}\in\sigma$$\hfill(1)$\\ Indeed: if $u\in 
\{x\}*(y\circ z)$, then $u\in x\circ t$ for some $t\in y\circ 
z\subseteq A_x$. Since $t\in A_x$, we have $(x, x\circ t)\in\sigma$. 
Since $u\in x\circ t$, we have $(x,u)\in\sigma$, so the property 
$(1)$ is satisfied. By $(1)$ and Lemma 77(1), we have ${\Big(}x\circ 
z,\, \{x\}*(y\circ z)*\{z\}{\Big)}\in\sigma$, that 
is,\begin{equation}\tag{$2$} {\Big(}x\circ z,\, (x\circ y)*(z\circ 
z{\Big)}\in\sigma\end{equation} On the other hand, since $(z, z\circ 
z)\in\sigma$, we have ${\Big(}(x\circ y)*\{z\},\, (x\circ y)*(z\circ 
z){\Big)}\in\sigma$.
In fact: Since $(z, z\circ z)\in\sigma$, by Lemma 2.10(2), we have
${\Big(}y\circ z, \{y\}*(z\circ z{\Big)}\in\sigma$.
Then, by Lemma 77(4), we have\begin{equation}\tag{$3$} 
{\Big(}\{x\}*(y\circ z), (x\circ y)*(z\circ 
z){\Big)}\in\sigma\end{equation}
By (3) and (2), we obtain ${\Big(}(x\circ y)*\{z\},\, x\circ 
z{\Big)}\in\sigma.$ Then, by $(1)$, we get $(x, x\circ z)\in\sigma$, 
and so $z\in A_x$. It remains to prove that $y\in A_x$.
Since $(x\circ z, x)\in\sigma$, by Lemma 77(1), we have 
${\Big(}(x\circ z)*\{y\}, x\circ y{\Big)}\in\sigma$. Since $(y\circ 
z, z\circ y)\in\sigma$, by Lemma 77(4), we have ${\Big(}\{x\}*(y\circ 
z), \{x\}*(z\circ y){\Big)}\in\sigma$. Thus we have 
${\Big(}\{x\}*(y\circ z), x\circ y{\Big)}\in\sigma$. Then, by $(1)$, 
we get $(x, x\circ y)\in\sigma$, and so $y\in A_x$.

Let now $y,z\in H$. Then we have $y\circ z\subseteq A_x$ or $(y\circ 
z)\cap A_x=\emptyset$. In fact:\\Let $y\circ z\nsubseteq A_x$ and 
$(y\circ z)\cap A_x\not=\emptyset$. Let $u\in y\circ z$ such that 
$u\notin A_x$, $v\in y\circ z$ and $v\in A_x$. Then we 
have\begin{equation}\tag{$1$}u\in y\circ z,\, (x, x\circ 
u)\notin\sigma,\, v\in y\circ z,\, (x, x\circ 
v)\in\sigma\end{equation}On the other 
hand,\begin{equation}\tag{$2$}(x, x\circ v)\in\sigma \mbox { and } 
v\in y\circ z \mbox { implies } {\Big(}x,\{x\}*(y\circ 
z){\Big)}\in\sigma\end{equation}Indeed: Let $a\in \{x\}*(y\circ z)$. 
Then $(x,a)\in\sigma$. Indeed: We have $a\in x\circ d$ for some $d\in 
y\circ z$. Since $(y\circ z, y\circ z)\in\sigma$, $v\in y\circ z$ and 
$d\in y\circ z$, we have $(v,d)\in\sigma$, then $(x\circ v, x\circ 
d)\in\sigma$. Since $(x, x\circ v)\in\sigma$ and $(x\circ v, x\circ 
d)\in\sigma$, we have $(x, x\circ d)\in\sigma$. Since $a\in x\circ 
d$, we have $(x,a)\in\sigma$.

In addition, since $\{u\}\subseteq y\circ z$, we have 
$\{x\}*\{u\}\subseteq \{x\}*(y\circ z)$. Since 
${\Big(}x,\{x\}*(y\circ z){\Big)}\in\sigma$, by Lemma 78, we have 
$(x, \{x\}*\{u\})\in\sigma$, that is $(x, x\circ u)\in\sigma$ which 
is impossible.

Since $A_x$ is a filter of $H$, by Proposition 79, we have 
$H\backslash A_x=\emptyset$ or $H\backslash A_x$ is a prime ideal of 
$H$. Then $H\backslash A_x=\emptyset$ or $H\backslash A_x$ is a 
proper prime ideal of $H$ (indeed, if  $H\backslash A_x=H$ then, 
since $A_x\subseteq H$, we have $A_x=\emptyset$ which is not 
possible).\\We consider the set$$\{H\backslash A_z \mid z\in H, \; 
H\backslash A_z \mbox { proper prime ideal of } H\}.$$We have $\sigma  
= \bigcap\limits_{z \in H} {{\sigma _{H\backslash {A_z}}}}$. In fact: 
Let $(x,y)\in\sigma$ and $z\in H$. Then $(x,y)\in\sigma_{H\backslash 
A_z}$. Indeed: Since $x\in H$, we have $x\in H\backslash A_z$ or 
$x\notin H\backslash A_z$. (a) If $x\notin H\backslash A_z$, then 
$x\in A_z$, so $(z, z\circ x)\in\sigma$. Since $(x,y)\in\sigma$, we 
have $(z\circ x, z\circ y)\in\sigma$. Then $(z, z\circ y)\in\sigma$, 
and $y\in A_z$, so $y\notin H\backslash A_z$. (b) If $x\in 
H\backslash A_z$. Then $y\in H\backslash A_z$. Indeed, if $y\notin 
H\backslash A_z$ then by (a), by symmetry, we have $x\notin 
H\backslash A_z$, which is impossible. Since both $x$ and $y$ belong 
to $H\backslash A_z$ or both do not belong to $H\backslash A_z$, we 
have $(x,y)\in\sigma_{H\backslash A_z}$.

Finally, let $(x,y)\in\sigma_{H\backslash A_z}$ for every $z\in H$. 
Then $(x,y)\in\sigma$. In fact:\\Since $x\in H$, we have 
$(x,y)\in\sigma_{H\backslash A_x}$. Since $x\in A_x$, we have 
$x\notin H\backslash A_x$. Since $(x,y)\in\sigma_{H\backslash A_x}$ 
and $x\notin H\backslash A_x$, we have $y\notin H\backslash A_x$, so 
$y\in A_x$. Since $y\in A_x$, we have $(x, x\circ y)\in\sigma$. 
Furthermore, since $(x,y)\in\sigma_{H\backslash A_y}$ and $y\in A_y$, 
that is $y\notin H\backslash A_y$, we have $x\notin H\backslash A_y$, 
so $x\in A_y$, and $(y, y\circ x)\in\sigma$. In addition, $\sigma$ is 
a semilattice congruence on $H$, so $(x\circ y, y\circ x)\in\sigma$. 
Since $(x, x\circ y)\in\sigma$,\, $(x\circ y, y\circ x)\in\sigma\,$ 
and $(y\circ x, y)\in\sigma$, we get $(x, y)\in\sigma$. $\hfill\Box$

The following proposition holds, and its proof is exactly as the 
proof of the Proposition in [10].\medskip

\noindent{\bf Proposition 82.} {\it Let H be an hypergroupoid and 
${\cal P}(H)$ the set of prime ideals of H. Then we have$${\cal N} = 
\bigcap\limits_{I \in {\cal P}(H)} {{\sigma _I}}.$$}{\bf Theorem 83.} 
{\it If H is an hypersemigroup, then the relation ${\cal N}$ is the 
least semilattice congruence on H.}\medskip

\noindent{\bf Proof.} Let $\sigma$ be a semilattice congruence on 
$H$. Then ${\cal N}\subseteq \sigma$. In fact: By Proposition 81, 
there exists a family $\cal A$ of proper prime ideals of $H$ such 
that $\sigma= \bigcap\limits_{I \in \cal A} {{\sigma _{\rm I}}}.$ By 
Proposition 82, ${\cal N} = \bigcap\limits_{I \in {\cal P}(H)} 
{{\sigma _I}}$, where ${\cal P}(H)$ is the set of prime ideals of 
$H$. On the other hand, $\bigcap\limits_{I \in \cal A} {{\sigma _{\rm 
I}}}\supseteq \bigcap\limits_{I \in {\cal P}(H)} {{\sigma _{\rm 
I}}}.$ Indeed: If $(x,y)\in\sigma_I$ for every prime ideal of $H$, 
then clearly $(x,y)\in\sigma_I$  for every proper prime ideal of $H$, 
and so for the elements of $\cal A$ as well. Hence we obtain ${\cal 
N}\subseteq \sigma$. $\hfill\Box$\medskip

\noindent{\bf Definition 84.} If $(H,\circ,\le)$ is an ordered 
hypersemigroup, and $F$ a filter of $(H,\circ)$, then $F$ is called a 
filter of $(H,\circ,\le)$ if $a\in F$ and $H\ni b\ge a$ implies $b\in 
F$. \medskip

\noindent{\bf Problem 85.} Using computer construct an ordered 
hypersemigroup $H$ for which the semilattice congruence $\cal N$ is 
not the least semilattice congruence on $H$.\medskip

\noindent{\bf Definition 86.} Let $(H,\circ)$ be an hypergroupoid 
endowed with a relation ``$\le$". We call it $\le$--hypergroupoid. A 
semilattice congruence $\sigma$ on $H$ is called {\it complete} if 
$(a,b)\in\le$ implies $(a,a\circ b)\in\sigma$, in the sense that if 
$u\in a\circ b$, then $(a,u)\in\sigma$. \medskip

\noindent{\bf Theorem 87.} {\it Let $(H,\circ)$ be an 
$\le$-hypergroupoid. Then the relation $\cal N$ is a complete 
semilattice congruence on $H$.}\medskip

\noindent{\bf Proof.} Let $a\le b$. We have to prove that $(a,a\circ 
b)\in{\cal N}$, that is, if $u\in a\circ b$, then $(a,u)\in\cal N$. 
Indeed: Since $N(a)\ni a\le b$, we have $b\in N(a)$. Since $a,b\in 
N(a)$, we have $a\circ b\subseteq N(a)$. Since $u\in a\circ b$, we 
have $u\in N(a)$, then $N(u)\subseteq N(a)$. On the other hand, since 
$u\in a\circ b$ and $u\in N(u)$, we have $u\in (a\circ b)\cap N(u)$. 
Since $(a\circ b)\not=\emptyset$, we have $a\circ b\subseteq N(u)$. 
Then $a\in N(u)$, and $N(a)\subseteq N(u)$. Hence we obtain 
$N(u)=N(a)$ and the proof is complete. $\hfill\Box$\smallskip
{\small

\end{document}